\documentclass[reqno]{amsproc}
\usepackage[latin1]{inputenc}
\usepackage{amssymb}
\usepackage{amsmath}
\usepackage{latexsym}%
\usepackage[frame,ps,matrix,arrow,curve,rotate,all,2cell,tips]{xy}
\usepackage{cite}
\usepackage{pgf,tikz}
\usepackage[normalem]{ulem}
\usetikzlibrary{matrix}
\usetikzlibrary{arrows}
\usepackage{amssymb}
\usepackage{amsfonts}
\usepackage{amscd}
\usepackage{xcolor}
\usepackage{hyperref}
\usepackage{amstext,amsmath,amssymb,amsfonts}
\newtheorem{theorem}{Theorem}

\newtheorem{corollary}[theorem]{Corollary}

\newtheorem{definition}[theorem]{Definition}
\newtheorem{example}[theorem]{Example}

\newtheorem{lemma}[theorem]{Lemma}

\newtheorem{proposition}[theorem]{Proposition}
\newtheorem{remark}[theorem]{Remark}

\newcommand{\K}{\mathbb {K}}

\newcommand{\A}{\mathcal{A}}

\newcommand{\G}{\mathcal{G}}

\newcommand{\B}{\mathcal{B}}
\newcommand{\bb}{\mathfrak{B}}

\newcommand{\h}{\mathcal{H}}

\newcommand{\beq}{\begin{eqnarray}}
\newcommand{\eeq}{\end{eqnarray}}
\newcommand{\beqs}{\begin{eqnarray*}}
\newcommand{\eeqs}{\end{eqnarray*}}
\newcommand{\bpro}{\begin{pro}}
\newcommand{\epro}{\end{pro}}
\newcommand{\blem}{\begin{lem}}
\newcommand{\elem}{\end{lem}}
\newcommand{\bdfn}{\begin{dfn}}
\newcommand{\edfn}{\end{dfn}}
\newcommand{\bcor}{\begin{cor}}
\newcommand{\ecor}{\end{cor}}
\newcommand{\bthm}{\begin{thm}}
\newcommand{\ethm}{\end{thm}}
\newcommand{\bex}{\begin{ex}}
\newcommand{\eex}{\end{ex}}
\newcommand{\brmk}{\begin{rmk}}
\newcommand{\ermk}{\end{rmk}}
\newcommand{\bpr}{\begin{pr}}
\newcommand{\epr}{\end{pr}}
\newcommand{\benum}{\begin{enumerate}} 
\newcommand{\eenum}{\end{enumerate}}
\newcommand{\bitem}{\begin{itemize}}
\newcommand{\eitem}{\end{itemize}}
\newcommand{\g}{\frak{g}}
\newcommand{\cqfd}{\hfill{\square}}
\chardef\bslash=`\\
\numberwithin{equation}{section}
\numberwithin{table}{section}
\numberwithin{theorem}{section}
\DeclareMathOperator{\id}{id}
\DeclareMathOperator{\ad}{ad}

\title{Hom-center-symmetric algebras and bialgebras}
\author{Mahouton Norbert Hounkonnou$^\ast$}
\address[$\ast, \dagger$]{University of Abomey-Calavi,
International Chair in Mathematical Physics and Applications,
ICMPA-UNESCO Chair, 072 BP 50, Cotonou, Rep. of Benin}
\email{$^\ast$ norbert.hounkonnou@cipma.uac.bj, with copy to hounkonnou@yahoo.fr}
\author{Mafoya Landry Dassoundo$^\dagger$}
\email{$^\dagger$ mafoya.dassoundo@cipma.uac.bj}
\begin{document}
\maketitle
\begin{abstract}
In this work,  the hom-center-symmetric algebras are constructed and discussed.  Their 
bimodules, dual bimodules and matched pairs are defined.
  The relation between the dual bimodules of  hom-center-symmetric algebras  and the matched pairs
of  hom-Lie algebras is established. 
Furthermore, the Manin triple of  hom-center-symmetric algebras is given. Finally, a
theorem linking the matched pairs of  hom-center-symmetric algebras, the
 hom-center-symmetric bialgebras and the matched
pairs of  sub-adjacent hom-Lie algebras
 is provided.
\end{abstract}


\noindent{
\today}
\tableofcontents
\section{Introduction}
The   Virasoro, Witt and  Lie algebra  deformations
generate  classes of nonassociative algebras \cite{Schafer} with   some interesting
algebraic identities (see \cite{Larsson_S1} and references therein). 
The Lie admissible algebras\cite{Santilli}  form  a class of nonassociative algebras with a product commutator defining a
 Lie algebra. 
The   hom-Lie algebras were first introduced in 2006 \cite{Hartwing_L_S},  using the $\sigma$-derivation of Virasoro and Witt
algebras. Some
$q$-deformations of Virasoro and Witt algebras also generate a  hom-Lie algebra structure. 
 The hom-Lie admissible algebras\cite{Larsson_S} constitute a class of 
nonassociative algebras with a product commutator giving a hom-Lie algebra.   

The algebraic properties of hom-coalgebras, hom-coassociative coalgebras,
and $G$ - hom - coalgebras, where $G$ is a subgroup  of the permutation group $S_3$ \cite{Makhlouf_S}
 generalizing  Lie-admissible coalgebras,  were investigated in \cite{Goze_R}.
  In these works, relevant definitions and properties of hom-Hopf algebras  generalizing  Hopf algebras, and
giving the module and comodule structures over hom-associative algebras and
hom-coassociative coalgebras,  were given.
Besides, the hom-Lie algebras  were studied in terms of  representation theory  \cite{Sheng1},
hom-Lie bialgebras\cite{Sheng_B}, hom-Lie 2-algebras\cite{Sheng_C}, enveloping algebras\cite{Yau}, and
Novikov algebras\cite{Yau2,Zhang_H_B }. 

Nearly hom-Lie algebras,  
left-symmetric algebras (LSA),  left-symmetric bialgebras (LSBA) were  analysed 
by analogy  to Lie bialgebras\cite{Drinfeld}.
C. Bai \cite{Bai} related left-symmetric algebras  with symplectic Lie algebras, classical Yang-Baxter equations
(CYBE), $\mathcal{O}$-operators, and para-K\"ahler Lie algebras. Similarly, hom-left-symmetric bialgebras
and para-K\"ahler hom-Lie algebras were equivalently introduced by Q. Sun and H. Li \cite{Sun_L}. 
The center-symmetric algebras were  investigated in \cite{Hounkonnou_D},  where  their  
Lie-admissibility was also established, and  their bimodules  were constructed. 
Furthermore, their matched pairs  were defined and 
linked to matched-pairs of Lie algebras and  associated  Manin triple.

 The hom-classical Yang-Baxter equation (HCYBE) was established from  a
homeomorphism derivation of the CYBE.  The Rota-Baxter hom-Lie admissible
algebra identities were also derived    and discussed in a series of works.   For more details, see  \cite{Yau1, Yau3, Yau4, Makhlouf_Y}  and references therein.

The   links between  hom-associative algebras (HAA), 
hom-associative bialgebras (HABA), hom-Lie algebras (HLA), hom-Lie bialgebras (HLBA),
hom-left-symmetric algebras (HLSA), and hom-left-symmetric bialgebras (HLSBA) can be represented by the following diagram: 
\begin{equation*}
\SelectTips{cm}{10}
\xymatrix{
 HAA \ar[r]^-{} \ar[d]_-{} & HLA \ar[d]^-{}& \ar[l]_-{} HLSA \ar[d]^-{}  \\
HABA\ar[r]_-{} & HLBA &\ar[l]_-{}  HLSBA.
  }
\end{equation*}

Our present work aims at investigating the  relations existing between a
 hom-center-symmetric algebra (HCSA), a  hom-center-symmetric bialgebra (HCSBA), a
 hom-Lie algebra (HLA) and a hom-Lie bialgebra (HLBA), as  illustrated by 
the following diagram:
 \begin{equation*}
     { 
     \xymatrix@C=40pt{
    HCSA \ar[r]^-{?} \ar[d]_-{?} & HLA \ar[d]^-{}  \\
   HCSBA\ar[r]_-{?} & HLBA }
     }
     \end{equation*} 
 We give  basic definitions  of hom-center-symmetric algebras and establish their  properties. 
 We show that the hom-center-symmetric algebras are  deformations (homomorphism deformations) of 
center-symmetric algebras. In addition, the bimodules and  matched pairs of hom-center-symmetric algebras are built and 
linked  to the bimodules and matched pairs of hom-Lie algebras, respectively. 
Besides, the Manin triple of hom-center-symmetric algebras is defined, 
discussed,  and analyzed with respect to the  Manin triple of hom-Lie algebras. 
Finally, a theorem yielding the link between the  hom-center-symmetric bialgebras, the
matched-pairs of hom-center-symmetric algebras, and the matched-pairs of hom-Lie algebras is given and proved. 
 \section{Preliminaries: main definitions  and  properties}
We now develop the basic  definitions and properties of hom-center-symmetric algebras.
\begin{definition}\label{dfn_hom_csa}
A hom-center-symmetric algebra is a triple $(\A, \mu, \alpha)$, where $\A$ is a vector space, 
$\mu:\A\otimes \A \rightarrow \A$ is a bilinear map, and $\alpha\in \mathfrak{gl}(\A)$ such that, for all $x, y, z\in \A$  
\beq
\alpha(\mu(x, y))=\mu(\alpha(x), \alpha(y)),
\eeq
\beq\label{eq_hom_def_relation}
(x, y, z)_{\alpha, \mu}=(z, y, x)_{\alpha, \mu},
\eeq
where $(x, y, z)_{\alpha, \mu}:=\mu(\mu(x, y), \alpha(z))-\mu(\alpha(x), \mu(y, z))$ is  called $\alpha$-associator associated to the map $\mu.$
\end{definition}
\begin{remark}
 The relation \eqref{eq_hom_def_relation}  is equivalent to 
\beqs
\mu\circ(\alpha\otimes\mu-\mu\otimes\alpha)=(\mu\circ\tau)\circ((\mu\circ\tau)\otimes\alpha -\alpha\otimes(\mu\otimes\tau)),
\eeqs 
which can be illustrated by the following commutative diagram:
\beqs
\label{diagram1}
\SelectTips{cm}{10}
\xymatrix{
\A \otimes \A \otimes \A \ar[rr]^-{\alpha\otimes \mu-\mu\otimes \alpha} \ar[d]_-{(\mu\circ \tau)\otimes\alpha-\alpha\otimes(\mu\circ\tau)}& &
\A \otimes \A \ar[d]^-{\mu }\\
\A\otimes\A \ar[rr]^-{\mu\circ\tau} & & \A
}
\eeqs 
{ where} $\tau$ is the exchange map on $\A\otimes \A$. 
\end{remark}
In  the sequel, as matter of notation simplification, we denote $(\A, \mu, \alpha)$ by $(\A, \alpha)$ or by $\A.$ 

 The  left $L$ and right $R$ representations of the  bilinear product  given on $\A$ are defined by
 \begin{eqnarray}
 L: \A & \longrightarrow & \mathfrak{gl}(\A)  \cr
  x  & \longmapsto & L_x:
  \begin{array}{ccc}
 \A &\longrightarrow & \A \cr 
  y & \longmapsto & \mu(x, y):=x \cdot y, 
   \end{array}
\end{eqnarray}
\begin{eqnarray}
    R: \A & \longrightarrow & \mathfrak{gl}(\A)  \cr
     x  & \longmapsto & R_x:
     \begin{array}{ccc}
    \A &\longrightarrow & \A \cr 
     y & \longmapsto & \mu(y, x):= y \cdot x.
      \end{array}
 \end{eqnarray}
 We infer the adjoint representation  $\ad: = L-R$ of the sub-adjacent Lie algebra $\G(\A)$ of a center-symmetric algebra $\A$  as follows:
  \begin{eqnarray}
    \ad: \A & \longrightarrow & \mathfrak{gl}(\A)  \cr
      x  & \longmapsto & \ad_x:
      \begin{array}{ccc}
     \A &\longrightarrow & \A \cr 
      y & \longmapsto & \mu(x, y)-\mu(y, x):=[x, y],
       \end{array}
    \end{eqnarray}
such that  
  $\displaystyle \forall \; x, y \in \A, \ad_x(y):=(L_x-R_x)(y).$
\begin{definition}
A hom-Lie algebra is a triple $(\g,[,]_{\g}, \alpha_{\g})$ consisting of 
 a {vector space} $\g$, 
an {algebra homomorphism} $\alpha_{\g}:\g\rightarrow \g$, and a skew-symmetric
bilinear map $[\cdot, \cdot]_{\g}:\g\otimes\g \rightarrow \g$ satisfying the following relation, 
for all $x, y, z\in \g:$
\beq\label{eq_hom_Jacobi_identity}
[\alpha_{\g}(x), [y,z]_{\g}]_{\g}+[\alpha_{\g}(y), [z,x]_{\g}]_{\g}+[\alpha_{\g}(z), [x,y]_{\g}]_{\g}=0,
\eeq 
called twisted Jacoby identity.
\end{definition}
\begin{definition}\label{dfn_rep_hom_Lie}
A representation of a hom-Lie algebra $(\g,[,]_{\g}, \alpha_{\g})$ on a vector space $V$ with respect
to $\psi\in \mathfrak{gl}(V)$ is a linear map $\rho:\g\rightarrow \mathrm{gl}(V)$ such that for all $x, y\in\g $,
the following relations:
\beq\label{eq_rep_Hom_Lie1}
\rho(\alpha_{\g}(x))\circ \psi=\psi\circ \rho(x),
\eeq
\beq\label{eq_rep_Hom_Lie2}
\rho([x,y]_{\g})\circ\psi =\rho(\alpha_{\g}(x))\circ \rho(y)-\rho(\alpha_{\g}(y))\circ\rho(x) 
\eeq
are satisfied.
\end{definition}
\begin{proposition}
Given a representation $(\rho, \psi, V)$ of a hom-Lie algebra $(\g, [,]_{\g}, \alpha_{\g})$, there is a hom-Lie algebra
on a semidirect sum vector space $\g\oplus V$ given by the following identities, {for all $x,y\in \g$ and  $u,v\in V$}:
\beq
(\alpha_{\g}\oplus\psi)(x+u)=\alpha_{\g}(x)+ \psi(u),
\eeq
\beq
[(x+u), (y+v)]_{\g\ltimes_{\rho} V}=[x, y]_{\g}+ \rho(x)v-\rho(y)u.
\eeq  
\end{proposition}
\textbf{Proof}.

According to the Definition~\ref{dfn_rep_hom_Lie}, it is clear that the product $[, ]_{\g\ltimes_{\rho}V}$ 
is bilinear on $\g\oplus V$ and $\forall x,y\in \g$ and $u, v\in V$,
\beqs
[(x+u), (y+v)]_{\g\ltimes_{\rho} V}&=&[x, y]_{\g}+ \rho(x)v-\rho(y)u 
=
-([y, x]_{\g}- \rho(x)v+\rho(y)u) \cr 
&=&-([y+v, x+u]_{\g\ltimes_{\rho} V}).
\eeqs
The bilinear product $[, ]_{\g\ltimes_{\rho}V}$ is skew symmetric.
Besides, we have 
\beqs
(\alpha_{\g}\oplus\psi)([(x+u), (y+v)]_{\g\ltimes_{\rho} V}
&=& (\alpha_{\g}\oplus\psi)([x, y]_{\g}+ \rho(x)v-\rho(y)u)\cr 
&=&
 \alpha_{\g}([x,y]_{\g})
 +\psi(\rho(x)v) -\psi(\rho(y)u)\cr 
&=& [\alpha_{\g}(x),\alpha_{\g}(y)]_{\g}+\rho(\alpha_{\g}(x))\circ \psi(v)-\rho(\alpha_{\g}(y))\circ \psi(u)\cr 
&=& [\alpha_{\g}(x)+\psi(u), \alpha_{\g}(y)+\psi(v)]_{\g\ltimes_{\rho} V}\cr
&=&[(\alpha_{\g}\oplus\psi)(x+u), (\alpha_{\g}\oplus\psi)(y+v)]_{\g\ltimes{\rho}V}.
\eeqs
In addition, we have  for all $x,y, z\in \g$ and $u, v, w\in V$:
\beqs
&&[(\alpha_{\g}\oplus\psi)(x+u),[y+v,z+w]_{\g\ltimes_{\rho}V}]_{\g\ltimes_{\rho}V}+
[(\alpha_{\g}\oplus\psi)(y+v),[z+w,x+u]_{\g\ltimes_{\rho}V}]_{\g\ltimes_{\rho}V}\cr 
&+&
[(\alpha_{\g}\oplus\psi)(z+w),[x+u,y+v]_{\g\ltimes_{\rho}V}]_{\g\ltimes_{\rho}V}=
[(\alpha_{\g}(x)+\psi(u),[y+v,z+w]_{\g\ltimes_{\rho}V}]_{\g\ltimes_{\rho}V}\cr 
&+&
[(\alpha_{\g}(y)+\psi(v),[z+w,x+u]_{\g\ltimes_{\rho}V}]_{\g\ltimes_{\rho}V}+
[(\alpha_{\g}(z)+\psi(w),[x+u,y+v]_{\g\ltimes_{\rho}V}]_{\g\ltimes_{\rho}V}\cr 
&=& [\alpha(x)+\psi(u), [y,z]_{\g}+\rho(y)w-\rho(z)v]_{\g\ltimes_{\rho}V}+
[\alpha(y)+\psi(v), [z,x]_{\g}+\rho(z)u-\rho(x)w]_{\g\ltimes_{\rho}V}\cr 
&+&[\alpha(z)+\psi(w), [x,y]_{\g}+\rho(x)v-\rho(y)u]_{\g\ltimes_{\rho}V}=
[\alpha(x)+[y,z]_{\g}]_{\g}+\rho(\alpha(x))(\rho(y)w-\rho(z)v)\cr 
&-&\rho([y,z]_{\g})\psi(u)+ 
[\alpha(y)+[z,x]_{\g}]_{\g}+\rho(\alpha(y))(\rho(z)u-\rho(x)w)-
\rho
([z,x]_{\g})\psi(v)\cr
&+&
[\alpha(z)+[x,y]_{\g}]_{\g}
+\rho(\alpha(z))(\rho(x)v-\rho(y)u)-
 \rho([x,y]_{\g})\psi(w)= 
[\alpha(x)+[y,z]_{\g}]_{\g}\cr 
&+&[\alpha(y)+[z,x]_{\g}]_{\g}+
[\alpha(z)+[x,y]_{\g}]_{\g}+
(\rho(\alpha(x))\rho(y)-\rho(\alpha(y))\rho(x))w+(\rho(\alpha(x))\rho(z)\cr 
&-&
\rho(\alpha(z))\rho(x))v
+(\rho(\alpha(y))\rho(z)-\rho(\alpha(z))\rho(y))u-\rho([y,z])\psi(u)-\rho([z,x])\psi(v)\cr 
&-&\rho([x,y])\psi(w)
=
\rho([x,y])\psi(w)
+\rho([y,z])\psi(u)-\rho([x,z])\psi(v)-\rho([y,z])\psi(u)\cr &-&
\rho([z,x])\psi(v)+\rho([x,u])\psi(w)
=0.
\eeqs 
Therefore, any representation  $(\rho, \psi, V)$ of a hom-Lie algebra $(\g, [,]_{\g},\alpha_{\g})$ 
leads to a  hom-Lie algebra on $\g\oplus V$.
$\cqfd$
\begin{proposition}\label{prop_Hom_Lie_Adm}
Let $(\A, \alpha)$ be a hom-center-symmetric algebra. The following conditions are satisfied:
\begin{itemize}
\item The commutator $[x,y]:=xy-yx$ defines a hom-Lie algebra $\G(\A)_{\alpha}:=(\A, [, ], \alpha)$. 
\item The couple $(\ad:= L-R, \alpha)$ gives a regular representation of the sub-adjacent hom-Lie algebra $\G(\A)_{\alpha}$.
\end{itemize}
\end{proposition}
\textbf{Proof}.

Consider the hom-center-symmetric algebra $(\A, \alpha, \mu)$.
It is clear from the definition that the commutator  $[,]_{\mu}$ associated to $\mu$ is bilinear and 
skew symmetric on $\A$. For all $x, y, z\in \A$, we have:
\beqs
&&[\alpha(x),[y,z]_{\mu}]_{\mu}+[\alpha(y),[z,x]_{\mu}]_{\mu}+[\alpha(z),[x,y]_{\mu}]_{\mu}
=[\alpha(x), \mu(y,z)-\mu(z,y)]_{\mu} \cr 
&+&[\alpha(y), \mu(z,x)-\mu(x,z)]_{\mu}+[\alpha(z), \mu(x,y)-\mu(y,x)]_{\mu} = [\alpha(x), \mu(y,z)]_{\mu}\cr 
&+&[\alpha(y), \mu(z,x)]_{\mu}
[\alpha(z), \mu(x,y)]_{\mu}- [\alpha(x), \mu(z,y)]_{\mu} 
 -[\alpha(y), \mu(x,z)]_{\mu}-[\alpha(z), \mu(y,x)]_{\mu} \cr 
&=&\mu(\alpha(x), \mu(y,z))+\mu(\alpha(y), \mu(z,x)) 
+ \mu(\alpha(z), \mu(x,y))-\mu(\alpha(x), \mu(z,y)) \cr 
&
-& \mu(\alpha(y), \mu(x,z))-\mu(\alpha(z), \mu(y,z))+\mu(\mu(y,z), \alpha(x))-\mu(\mu(z,x), \alpha(y)) \cr 
&-& \mu(\mu(x,y), \alpha(z))+\mu(\mu(z,y), \alpha(x))+ \mu(\mu(x,z), \alpha(y))+\mu(\mu(y,x), \alpha(z)) \cr
&
=&\{\mu(\alpha(x), \mu(y,z))-\mu(\mu(x,y), \alpha(z)) \}
+\{\mu(\alpha(y), \mu(z,x))-\mu(\mu(y,z), \alpha(x)) \} \cr 
&+&\{\mu(\alpha(z), \mu(x,y))-\mu(\mu(z,x), \alpha(y)) \}
-\{\mu(\alpha(x), \mu(z,y))-\mu(\mu(x,z), \alpha(y)) \}\cr 
& -& \{\mu(\alpha(y), \mu(x,z))
-\mu(\mu(y,x), \alpha(z)) \}-\{\mu(\alpha(z), \mu(y,x))-\mu(\mu(z,y), \alpha(x)) \}\cr 
&=&
\{\mu(\alpha(x), \mu(y,z))-\mu(\mu(x,y), \alpha(z)) \}-\{\mu(\alpha(z), \mu(y,x))-\mu(\mu(z,y), \alpha(x)) \}\cr 
&+&\{\mu(\alpha(y), \mu(z,x))-\mu(\mu(y,z), \alpha(x)) \}-\{\mu(\alpha(x), \mu(z,y))-\mu(\mu(x,z), \alpha(y)) \}\cr 
&+&\{\mu(\alpha(z), \mu(x,y))-\mu(\mu(z,x), \alpha(y)) \}-\{\mu(\alpha(y), \mu(x,z))-\mu(\mu(y,x), \alpha(z)) \}\cr 
& 
=
& [\alpha(x),[y,z]_{\mu}]_{\mu}+[\alpha(y),[z,x]_{\mu}]_{\mu}+[\alpha(z),[x,y]_{\mu}]_{\mu}
= 0.
\eeqs
Therefore, for all $x, y, z\in\A$, the following equality holds
\beqs
[\alpha(x),[y,z]_{\mu}]_{\mu}+[\alpha(y),[z,x]_{\mu}]_{\mu}+[\alpha(z),[x,y]_{\mu}]_{\mu}=0,
\eeqs
and  shows that $\G(\A)_{\alpha}$ is a hom-Lie algebra.
$\cqfd$
\begin{remark}
The hom-Lie algebra $\G(\A)_{\alpha}:=(\A, [, ], \alpha)$ given in the Proposition~\ref{prop_Hom_Lie_Adm} is called the sub-adjacent hom-Lie algebra of 
$(\A, \alpha),$ and  $(\A, \alpha)$ is said to be  the compatible 
hom-center-symmetric algebra structure on the hom-Lie algebra $\G(\A)_{\alpha}$.
\end{remark}
 \section{Bimodule and matched pair of hom-center-symmetric algebras}
\begin{definition}
Consider a hom-center-symmetric algebra $(\A, \alpha),$  a vector space $V,$   two linear maps  
$l, r: \A \rightarrow \mathfrak{gl}(V),$ and $\varphi \in \mathfrak{gl}(V)$. The quadruple $(l,r,V, \varphi)$
is called bimodule of the hom-center-symmetric algebra $(\A, \alpha)$ if for all $x, y \in \A$ and $v \in V$,
\beq\label{bimod_hom_center1}
\varphi\circ l_x=l_{\alpha(x)}\circ\varphi,\quad\quad
 \varphi\circ r_x=r_{\alpha(x)}\circ\varphi,
\eeq
\beq\label{bimod_hom_center2}
l_{\alpha(x)}\circ l_y-l_{xy}\circ \varphi=r_{yx}\circ \varphi-r_{\alpha(x)}\circ r_y,
\eeq
\beq\label{bimod_hom_center3}
l_{\alpha(x)}\circ r_y-r_{\alpha(y)}\circ l_x=l_{\alpha(y)}\circ r_x-r_{\alpha(x)}\circ l_y.
\eeq
\end{definition} 
In the particular case when $\alpha=\id,$  the construction of the  center-symmetric algebra structure on the  
semi-direct vector space $\A\oplus V$ is given in \cite{Hounkonnou_D}. By analogy to this work, we have:
\begin{proposition}
Let $(\A, \alpha)$ be a hom-center-symmetric algebra, $V$ a vector space,
 $l, r:\A \rightarrow\mathfrak{gl}(V)$ two linear maps,  and $\varphi\in\mathfrak{gl}(V)$. 
Then the  quadruple $( \varphi, l, r, V)$ is a bimodule of the hom-center-symmetric algebra 
$(\A, \alpha)$ if and only if $(\A\oplus V, \ast, \alpha\oplus\varphi)$ is a 
hom-center-symmetric algebra, where $\ast $ and $\alpha\oplus \varphi$ are defined as follows:
for all $x, y\in \A$, $u, v\in V$,
\beqs
(x+u)\ast(y+v)=xy+l(x)v+r(y)u,
\eeqs
\beqs
(\alpha\oplus\varphi)(x+u)=\alpha(x)+\varphi(u).
\eeqs
\end{proposition}
\textbf{Proof}.

For all $x, y, z\in \A$ and  $u, v, w \in V$, using the relations \eqref{bimod_hom_center1}, we have:
\beqs
(\alpha\oplus \varphi)((x+u)\ast (y+v))&=&(\alpha\oplus \varphi)(xy+l_xv+r_yu)
=\alpha(xy)+\varphi(l_xv)+\varphi(r_yu)\cr 
&=&
\alpha(x)\alpha(y) 
+(\varphi\circ l_x)v+(\varphi\circ r_y)u\cr 
&=&\alpha(x)\alpha(y)+(l\circ\alpha(x)) \varphi (v)+(r\alpha(y))\ \varphi (u)
\cr &=& (\alpha(x)+\varphi(u))\ast (\alpha(y)+\varphi(v)).
\eeqs
The associator associated to the bilinear product $\ast$ gives:
\beqs
&&(x+u, y+v, z+w)
=((x+u)\ast (y+v))\ast (\alpha\oplus \varphi)(z+w)
\cr 
&-&
(\alpha\oplus \varphi)(x+u)\ast ((y+v)\ast (z+w))=(xy+l_xv+r_yu)\ast(\alpha(z)+\varphi(w))
\cr
&-&(\alpha(x)+\varphi(u))\ast (yz+l_yw+r_zv)= (xy)\alpha(z)+l_{xy}\varphi(w)
+(r\circ\alpha(z))(l_xv+r_yu)-\alpha(x)(yz)\cr 
&-&r_{yz}\varphi(u)
-((l\circ\alpha(x))(l_yw+r_zv))=((xy)\alpha(z)-\alpha(x)(yz))\cr 
&+&(l_{xy}\circ\varphi-(l\circ\alpha(x))\circ l_y)w
+((r\circ\alpha(z))\circ l_x-l\circ\alpha(x)r_z)v+(r_\circ\alpha(z)r_y-r_{yz})u
= (x, y,z)_{\alpha}\cr 
&+&(l_{xy}\circ\varphi-(l\circ\alpha(x))\circ l_y)w+
((r\circ \alpha(z))\circ l_x-(l\circ\alpha(x))\circ r_z)v+
((r\circ\alpha(z))\circ r_y-r_{yz})u.
\eeqs
We also have:
\beqs
(z+w, y+v, x+u)_{_{\alpha\oplus\varphi}}
&=& (z, y,x)_{\alpha}+(l_{zy}\circ\varphi-(l\circ\alpha(z))\circ l_y)u+  
((r\circ \alpha(x))\circ l_z\cr&&
-(l\circ\alpha(z))\circ r_x)v  +
((r\circ\alpha(x))\circ r_y-r_{yx})w.
\eeqs
By using the relations \eqref{bimod_hom_center2} and \eqref{bimod_hom_center3},  we have: 
\beqs
&&(x+u, y+v, z+w)_{_{\alpha\oplus\varphi}}-(z+w, y+v, x+u)_{_{\alpha\oplus\varphi}}
=
(x, y,z)_{\alpha}-(z, y,x)_{\alpha}
\cr &+&
((l_{zy}\circ\varphi-(l\circ\alpha(z))\circ l_y)
-((r\circ\alpha(z))\circ r_y-r_{yz}))u 
+(((r\circ \alpha(z))\circ l_x-(l\circ\alpha(x))\circ r_z)\cr 
&-&((r\circ \alpha(z))\circ l_x-(l\circ\alpha(x))\circ r_z))v
+((l_{xy}\circ\varphi-(l\circ\alpha(x))\circ l_y)
-((r\circ\alpha(x))\circ r_y-r_{yx}))w \cr
&=& 0.
\eeqs
Therefore, $( \varphi, l, r, V)$ is a bimodule of the hom-center-symmetric algebra $(\A, \alpha)$ if and only
if the equations \eqref{bimod_hom_center1}, \eqref{bimod_hom_center2}, and \eqref{bimod_hom_center3} are satisfied.
$\cqfd$

We denote the hom-center-symmetric algebra $(\A\oplus V, \ast, \alpha\oplus\varphi)$  
by $\A\ltimes_{l,r}^{\alpha, \varphi}V$ or simply by $\A\ltimes V$.
\begin{example}
Let $(\A, \alpha)$ be a hom-center-symmetric algebra. For all $x, y, z \in \A$, we have:
\beqs
(x,y,z)_{\alpha}=(z,y,x)_{\alpha} 
&\Longleftrightarrow&
(xy)\alpha(z)-\alpha(x)(yz)=(zy)\alpha(x)-\alpha(z)(yx)\cr 
&\Longleftrightarrow& (L_{xy}\circ\alpha-L_{\alpha(x)}\circ L_{y})(z)=(R_{\alpha(x)}\circ R_y-R_{yx}\circ \alpha)(z)\cr 
&\Longleftrightarrow& ((L_{xy}\circ\alpha-L_{\alpha(x)}\circ L_{y})-(R_{\alpha(x)}\circ R_y-R_{yx}\circ \alpha))(z)=0.
\eeqs
\beq\label{eq_hom_center_example1}
(x,y,z)_{\alpha}=(z,y,x)_{\alpha} 
&\Longleftrightarrow& ((L_{xy}\circ\alpha-L_{\alpha(x)}\circ L_{y})-(R_{\alpha(x)}\circ R_y-R_{yx}\circ \alpha))(z)=0.
\eeq
On the other hand, we also have
\beqs
(x,y,z)_{\alpha}=(z,y,x)_{\alpha} 
&\Longleftrightarrow& (xy)\alpha(z)-\alpha(x)(yz)=(zy)\alpha(x)-\alpha(z)(yx)\cr 
&\Longleftrightarrow&((R_{\alpha(z)}\circ L_x)-(L_{\alpha(x)}\circ R_z)-
((R_{\alpha(x)}\circ L_z)-(L_{\alpha(z)})\circ R_x)(y)=0. 
\eeqs
\beq\label{eq_hom_center_example2}
(x,y,z)_{\alpha}=(z,y,x)_{\alpha} 
\Longleftrightarrow((R_{\alpha(z)}\circ L_x)-(L_{\alpha(x)}\circ R_z)-
((R_{\alpha(x)}\circ L_z)-(L_{\alpha(z)})\circ R_x)=0. 
\eeq
From   the relations \eqref{eq_hom_center_example1} and \eqref{eq_hom_center_example2}, where the 
right hand sides translate that $(\alpha, L, R)$ satisfies the relations \eqref{bimod_hom_center2} and \eqref{bimod_hom_center3},
 we deduce that 
$( \alpha, L,R, \A)$ is a bimodule of $(\A, \alpha)$.  
\end{example}

\begin{definition}
A homomorphism of two hom-center-symmetric algebras $(\A_1,\mu_1, \alpha_1)$ and $(\A_2,\mu_2, \alpha_2)$ is a linear
map $f: \A_1\rightarrow \A_2$ such that:
\beq
f\circ\mu_1=\mu_2\circ(f\otimes f),
\eeq
\beq
f\circ\alpha_1=\alpha_2\circ f,
\eeq
as illustrated, respectively,  by the following commutative diagrams:
\beqs
\label{diagram2}
\SelectTips{cm}{10}
\xymatrix{
\A_1\otimes \A_1 \ar[rr]^-{\mu_1} \ar[d]_-{f\otimes f}& & \A_1 \ar[d]^-{f}\\
\A_2\otimes \A_2\ar[rr]^-{\mu_2} & & \A_2}
\quad \quad \quad \quad
\SelectTips{cm}{10}
\xymatrix{
\A_1 \ar[rr]^-{f} \ar[d]_-{\alpha_1}& & \A_2 \ar[d]^-{\alpha_2}\\
\A_1\ar[rr]^-{f} & & \A_2}
\eeqs 
\end{definition}
\begin{lemma}\label{lemma_hom_1}
Let $( \varphi, l, r, V)$ be a bimodule of a hom-center-symmetric algebra $(\A, \alpha)$. Then we have:
\begin{itemize}
\item The couple $(\rho=l-r, \varphi)$ is a representation of the sub-adjacent hom-Lie algebra $\G(\A)$ associated to $(\A, \alpha);$
\item For any representation $(\rho, \psi),$ with $\rho:\G(\A)\rightarrow \mathfrak{gl}(\G(\A)),$ of { the underlying} 
hom-Lie algebra $(\G(\A), \alpha)$ { of the hom-center-symmetric algebra $(\A, \alpha)$}, the triple $(\rho, 0, \psi)$ is a bimodule of $(\A,\psi);$
\item The hom-center-symmetric algebras $\A\ltimes_{l, r}^{\alpha, \varphi} V$ and 
$\A\ltimes_{l-r, 0}^{\alpha, \varphi}V$ given by the bimodules $(l,r,\varphi)$ and $(l-r, 0, \varphi)$, respectively,
have the same sub-adjacent hom-Lie algebra given by the semidirect sum $\G(\A)\ltimes_{l-r}^{\alpha, \varphi} V$ of the hom-Lie algebra $(\G(\A), \alpha)$ and its representation $(l-r, \alpha, V)$ as:
$
[x+u, y+v]=[x,y]+(l-r)(x)v-(l-r)(y)u, \quad \forall x, y \in \A, u, v\in V.
$ 
\end{itemize}
\end{lemma}
\textbf{Proof}.

Let $( \varphi, l, r, V)$ be a bimodule of a hom-center-symmetric algebra $(\A, \alpha)$.
We have for all $x, y, z\in \A$ and  $u, v,w\in V$,
\begin{itemize}
\item 
$(l-r)(\alpha(x))\circ\varphi=l(\alpha(x))\circ\varphi-r(\alpha(x))\circ\varphi=\varphi\circ l_x-\varphi\circ r_x=\varphi\circ (l-r)(x)$.
\beqs
&&(l-r)([x,y])\circ\varphi=(l-r)(xy)\circ\varphi-(l-r)(yx)\circ\varphi 
= l(xy)\circ\varphi-r(xy)\circ\varphi-l(yx)\circ\varphi\cr 
&+& r(yx)\circ\varphi= (l(xy)+r(yx))\varphi-(l(yx)+r(xy))\varphi=(l\circ\alpha(x))\circ l(y)+(r\circ\alpha(x))\circ r(y)\cr 
&-&(l\circ\alpha(y))\circ l(x)-(r\circ\alpha(y))\circ r(x)= (l\circ\alpha(x))\circ l(y)+(r\circ\alpha(x))\circ r(y)
-(l\circ\alpha(y))\circ l(x)\cr 
&-&(r\circ\alpha(y))\circ r(x)-(l\circ \alpha(x))\circ r(y)+(r\circ \alpha(y))\circ l(x)+(l\circ \alpha(y))\circ r(x)-(r\circ \alpha(x))\circ l(y) 
\cr 
&
=& ((l-r)\circ\alpha(x))\circ(l-r)(y)-((l-r)\circ\alpha(y)\circ(l-r)(x).
\eeqs
Therefore, $(l-r, \varphi)$ is a representation of the sub-adjacent hom-Lie algebra $\G(\A)_{\alpha}$ of the hom-center-symmetric algebra $(\A,\alpha)$.
\item 
Using the fact that  
$(\rho, \psi )$ is a representation of $\G(\A)_{\alpha}$ and setting
$r=0$, then the relations \eqref{bimod_hom_center3} and \eqref{bimod_hom_center1} are satisfied,  and, in addition, we have,  for all $x,y\in \A$,
\beqs
\rho([x,y])\circ \psi=\rho(\alpha(x))\circ\rho(y)-\rho(\alpha(y))\circ\rho(x)
\eeqs
which implies 
\beqs
\rho(xy)\circ\psi-\rho(\alpha(x))\circ\rho(y)=\rho(yx)\circ\psi-\rho(\alpha(y))\circ\rho(x).
\eeqs
This is exactly the equation\eqref{bimod_hom_center2}.
\item The commutator of the bilinear product defined as  $(x+u)\ast(y+v):=xy+l_xv+r_yu$ is given by
$(x+u)\ast( y+v)-(y+v)\ast (x+u) =[x,y]+(l-r)(x)v-(l-r)(y)u=[x+u, y+v]$.
\end{itemize}
$\cqfd$
\begin{definition} 
Let $(\A, \alpha)$ be a hom-center-symmetric algebra.
The dual maps $l^{*}, r^{*}$  of the linear maps $l,r:\A\rightarrow \mathfrak{gl}(V)$ are defined, respectively, as: 
\begin{equation}
\begin{aligned}
\label{dual_linear_1}
 l^*: \A & \longrightarrow  \mathfrak{gl}(V^*)  \\
  x  & \longmapsto  l^*_x:
      \begin{array}{llll}
 V^* &\longrightarrow & V^* \\ 
  u^* & \longmapsto & l^*_x u^*: 
      \begin{array}{llll}
V  &\longrightarrow&  \K \\
v  &\longmapsto& \left< l^{*}_xu^{*}, v \right>
 :=   \left<  u^{*}, l_x v\right>, 
      \end{array}
     \end{array}
\end{aligned}
\end{equation}
\begin{equation}
\begin{aligned}\label{dual_linear_2}
 r^*: \A & \longrightarrow  \mathfrak{gl}(V^*)  \\
  x  & \longmapsto  r^*_x:
      \begin{array}{llll}
 V^* &\longrightarrow & V^* \\ 
  u^* & \longmapsto & r^*_x u^*: 
      \begin{array}{lllll}
V  &\longrightarrow&  \K \\
v  &\longmapsto& \left< r^{*}_xu^{*}, v \right>
 :=  \left<  u^{*}, r_x v\right>, 
      \end{array}
     \end{array}
\end{aligned}
\end{equation}
for all $x \in \A, u^{*} \in V^{*}, v \in V.$ 
\end{definition}
\begin{proposition}
Let $(\A, \alpha)$ be a hom-center-symmetric algebra, and let $l, r : \A \rightarrow \mathfrak{gl}(V )$
and $\varphi: V\rightarrow V$, where $V$ is a finite dimensional vector space, be three linear maps. Then, the following
conditions are equivalent:
\begin{itemize}
\item $(l, r, \varphi, V )$ is a bimodule of $\A$.
\item $(r^* , l^* , \varphi^*, V^* )$ is a bimodule of $\A$ and  $\alpha^2=\id$.
\end{itemize} 
\end{proposition}
\textbf{Proof}.

Consider a hom-center-symmetric algebra $(\A, \alpha)$ and two linear maps
 $l, r: \A\rightarrow \mathfrak{gl}(V),$  with their dual maps $l^*, r^*$ satisfying  the relations
\eqref{dual_linear_1} and \eqref{dual_linear_2}.
\begin{itemize}
\item  $\Longrightarrow)$
 At first, let's suppose that $(l,r, V)$ is a bimodule of the hom-center-symmetric algebra $(\A, \alpha)$.
For all $x, y\in \A$,  $v\in V$, and  $u^*\in V^*$, we have:
\beqs
\left<(\varphi^*\circ r_x^*)u^*, v  \right>&=& \left<r_x^*u^*, \varphi(v)  \right>=\left<u^*, r_x\varphi(v)  \right>\crcr&=&
\left<u^*, \varphi\circ r_{\alpha(x)}v  \right>=\left<\varphi^*(u^*), r_{\alpha(x)}v  \right>
\eeqs
\beqs
\left<(\varphi^*\circ l_x^*)u^*, v  \right>&=& \left<l_x^*u^*, \varphi(v)  \right>=\left<u^*, l_x\varphi(v)  \right>
\cr
&=&
\left<u^*, \varphi\circ l_{\alpha(x)}v  \right>=\left<\varphi^*(u^*), l_{\alpha(x)}v  \right> 
\eeqs
Thus,
\beqs
\varphi^*\circ l_x^*=l^*_{\alpha(x)}\circ \varphi^*, \quad\quad \varphi^*\circ r_x^*=r^*_{\alpha(x)}\circ \varphi^*.
\eeqs
Besides,
\beqs
\left<(r^*_{\alpha(x)}r^*_y-r_{xy}^*\circ\varphi^*)u^*, v\right>&=&\left<u^*, (r_yr_{\alpha(x)}-\varphi\circ r_{xy})v\right>= \left<u^*, (r_{\alpha^2(y)}r_{\alpha(x)}-r_{\alpha(xy)}\circ\varphi)v\right>\cr
&=& \left<u^*, (r_{\alpha^2(y)}r_{\alpha(x)}-r_{\alpha(x)\alpha(y)}\circ\varphi)v\right>\cr 
&=& \left< u^*,  (l_{\alpha(yx)}\circ \varphi-l_{\alpha^2(y)}l_{\alpha(x)})v   \right>= \left<u^*, (\varphi\circ l_{yx}-l_{y}l_{\alpha(x)})v\right>\cr 
\left<(r^*_{\alpha(x)}r^*_y-r_{xy}^*\circ\varphi^*)u^*, v\right>
&=& \left<(l^*_{yx}\circ \varphi^*-l^*_{\alpha(x)}l^*_{y})u^*, v\right>.
\eeqs
It follows that 
\beq\label{eq_dual_homm_1}
r^*_{\alpha(x)}r^*_y-r_{xy}^*\circ\varphi^*=l^*_{yx}\circ\varphi^*-l^*_{\alpha(x)}l^*_{y}.
\eeq
Setting  $\alpha^2=\id_{\A}$ leads to
\beqs
\left<(r^*_{\alpha(x)}l^*_y-l^*_{\alpha(y)}r^*_x)u^* , v\right>&=&\left<u^*,(l_yr_{\alpha(x)}-r_xl_{\alpha(y)}) v\right>=\left<u^* , (l_{\alpha^2(y)}r_{\alpha(x)}-r_{\alpha^2(x)}l_{\alpha(y)}) v\right>\cr 
&=&\left<u^* , (l_{\alpha^2(x)}r_{\alpha(y)}-r_{\alpha^2(y)}l_{\alpha(x)}) v\right>
=\left<u^* , (l_xr_{\alpha(y)}-r_{y}l_{\alpha(x)}) v\right>\cr  
&=&\left<(r^*_{\alpha(y)}l^*_x-l^*_{\alpha(x)}r^*_{y})u^*,v\right>.
\eeqs
Then, using $\alpha^2=\id_{\A}$,  the relation
\beq\label{eq_dual_homm_2}
r^*_{\alpha(x)}l^*_y-l^*_{\alpha(y)}r^*_x=r^*_{\alpha(y)}l^*_x-l^*_{\alpha(y)}r^*_{x}
\eeq
is satisfied. 

Finally, from the relations \eqref{eq_dual_homm_1} and \eqref{eq_dual_homm_2}, we conclude that the triple
$(r^*, l^*, V^*)$ is a bimodule of the hom-center-symmetric algebra $(\A, \alpha)$.
\item  $\Longleftarrow) $
Conversely, suppose that $\alpha^2=\id_{\A}$ and $(r^*, l^*, V^*)$ is a bimodule of $\A$.
Then, it is straightforward to   check that the triple  $(l,r, V)$ is a bimodule of $\A.$ 
\end{itemize}
Therefore, it is true that the triple $(l, r, V)$ is a bimodule of the hom-center-symmetric algebra $(\A, \alpha)$ if, 
and only if,  the  triple $(r^*, l^*, V^*),$ where $r^*, l^*$ are  given by the relations 
\eqref{dual_linear_1} and \eqref{dual_linear_2}, respectively, is  a bimodule of the   
hom-center-symmetric algebra $(\A, \alpha)$ with $\alpha^2=\id_{\A}$.
$\cqfd$
\begin{theorem}\label{thm_matched_hom_center}
Let $(\A, \cdot, \alpha_{\A})$ and $(\B, \circ, \alpha_{\B})$ be two hom-center symmetric algebras.
Suppose  there are linear maps $l_{\A}, r_{\A}:\A\rightarrow \mathfrak{gl}(\B)$ 
and $l_{\B}, r_{\B}:\B\rightarrow \mathfrak{gl}(\A)$ such that $(l_{\A}, r_{\A}, \alpha_{\B})$ and 
$(l_{\B}, r_{\B}, \alpha_{\A})$ are  bimodules of the hom-center-symmetric algebra $\A$ and $\B$, respectively,
satisfying the following conditions for all $x, y\in \A$ and  $a, b \in \B:$
\beq\label{eq_matched_hom_center_1}
(l_{\B}(a)x)\cdot (\alpha_{\A}(y))+l_{\B}(r_{\A}(x)a)(\alpha_{\A}(y))-l_{\B}(\alpha_{\B}(a))(x\cdot y)
-r_{\B}(\alpha_{\B}(a))(y\cdot x)\cr 
+(\alpha_{\A}(y))\cdot (r_{\B}(a)x)+r_{\B}(l_{\A}(x)a)(\alpha_{\A}(y))=0,
\eeq
\beq\label{eq_matched_hom_center_2}
(r_{\B}(a)x)\cdot(\alpha_{\A}(y))+l_{\B}(l_{\A}(x)a)(\alpha_{\A}(y))-(\alpha_{\A}(x))\cdot (l_{\B}(a)y)
-r_{\B}(r_{\A}(y)a)(\alpha_{\A}(x))\cr 
-(r_{\B}(a)y)\cdot(\alpha_{\A}(x))-l_{\B}(l_{\A}(y)a)(\alpha_{\A}(x))+(\alpha_{\A}(y))\cdot (l_{\B}(a)x)
+r_{\B}(r_{\A}(x)a)(\alpha_{\A}(y))=0, 
\eeq
\beq\label{eq_matched_hom_center_3}
(l_{\B}(x)a)\circ (\alpha_{\B}(b))+l_{\A}(r_{\B}(a)x)(\alpha_{\B}(b))-l_{\A}(\alpha_{\A}(x))(a\circ b)
-r_{\A}(\alpha_{\A}(x))(b\circ a)\cr 
+(\alpha_{\B}(b))\circ (r_{\A}(x)a)+r_{\A}(l_{\B}(a)x)(\alpha_{\B}(b))=0,
\eeq
\beq\label{eq_matched_hom_center_4}
(r_{\A}(x)a)\circ(\alpha_{\B}(b))+l_{\A}(l_{\B}(a)x)(\alpha_{\B}(b))-(\alpha_{\B}(a))\circ (l_{\A}(x)b)
-r_{\A}(r_{\B}(b)x)(\alpha_{\B}(a))\cr 
-(r_{\A}(x)b)\circ(\alpha_{\B}(a))-l_{\A}(l_{\B}(b)x)(\alpha_{\B}(a))+(\alpha_{\B}(b))\circ (l_{\A}(x)a)
+r_{\A}(r_{\B}(a)x)(\alpha_{\B}(b))=0.
\eeq
Besides, there exists  a hom-center-symmetric algebra structure on the vector space $\A\oplus \B$ given by:
\beq\label{eq_hom_csa}
(x+a)\ast(y+b)=(x\cdot y+l_{\B}(a)y+r_{\B}(b)x)+(a\circ b+l_{\A}(x)b+r_{\A}(y)a).
\eeq
\end{theorem}
\textbf{Proof}.

From its definition given in \eqref{eq_hom_csa}, the product $\ast$    is bilinear on $\A\oplus \B,$ and we have
for all $x,y\in \A$ and  $a, b\in \B:$
\beqs
&&(\alpha_{\A}+\alpha_{\B})((x+a)\ast(y+b))= \alpha_{\A}(x\cdot y+l_{\B}(a)y+r_{\B}(b)x)+
\alpha_{\B}(a\circ b+l_{\A}(x)b+r_{\A}(y)a)\cr 
&=&  \alpha_{\A}(x\cdot y)+ \alpha_{\A}(l_{\B}(a)y)+ \alpha_{\A}(r_{\B}(b)x)+
\alpha_{\B}(a\circ b) 
+\alpha_{\B}(l_{\A}(x)b)+ \alpha_{\B}(r_{\A}(y)a) \cr
&=& 
( \alpha_{\A}(x)\cdot \alpha_{\A}(y) +l_{\B}( \alpha_{\B}(a))  \alpha_{\A}(y)+r_{\B}( \alpha_{\B}(b)) \alpha_{\A}(x)) 
+( \alpha_{\B}(a)\circ \alpha_{\B}(a)+l_{\A}( \alpha_{\A}(x))  \alpha_{\B}(b)\cr
&+& r_{\A}( \alpha_{\A}(y)) \alpha_{\B}(a))
=(( \alpha_{\A}(x)+( \alpha_{\B}(a))\ast (( \alpha_{\A}(y)+( \alpha_{\B}(b)).
\eeqs
Using the trilinearity of the associator associated to the bilinear product $\ast$ on the vector space 
$\A\oplus \B$, we have, for all $x, y, z\in \A$ and  $a, b,c\in \B:$
\beq\label{eq_hom_asso_csa}
(x+a,y+b,z+c)_{\alpha_{\A} \oplus \alpha_{\B} }=(x,y,z)_{\alpha_{\A}}+
(x,y,c)_{\alpha_{\A} \oplus \alpha_{\B} }+
(x,b, z)_{\alpha_{\A} \oplus \alpha_{\B} }+(x,b, c)_{\alpha_{\A} \oplus \alpha_{\B} }\nonumber\\
+(a, y,z)_{\alpha_{\A} \oplus \alpha_{\B} }+ 
(a, y,c)_{\alpha_{\A} \oplus \alpha_{\B} }
+(a,b, z)_{\alpha_{\A} \oplus \alpha_{\B} }+ (a,b,c)_{\alpha_{\B}}.
\eeq
From the Definition~\ref{dfn_hom_csa}, the  hom-center-symmetric algebra structures, defined on 
$(\A, \alpha_{\A})$ and $(\B, \alpha_{\B}),$
are 
 $(x,y,z,)_{\alpha_{\A}}=(z,y,x)_{\alpha_{\A}}$, and 
$(a,b,c)_{\alpha_{\B}}=(c,b,a)_{\alpha_{\B}},$ respectively.

Before showing that the remaining associators given in  the equation \eqref{eq_hom_asso_csa} are center symmetric, 
it is useful to remark that the trilinear identities
$(x, b, c)_{\alpha_{\A} \oplus \alpha_{\B} }$, $(a, b, z)_{\alpha_{\A} \oplus \alpha_{\B} }$
, $(a, y, z)_{\alpha_{\A} \oplus \alpha_{\B} }$, and  $(x,y,c)_{\alpha_{\A} \oplus \alpha_{\B} }$
play a  symmetric role, what means that showing   $(x, b, c)_{\alpha_{\A} \oplus \alpha_{\B} }$ is 
center-symmetric   is equivalent to showing that  $(a, b, z)_{\alpha_{\A} \oplus \alpha_{\B} }$ is 
center-symmetric  too. Similarly, proving that $(a, y, z)_{\alpha_{\A} \oplus \alpha_{\B} }$ is
center-symmetric  is equivalent to proving that $(x,y,c)_{\alpha_{\A} \oplus \alpha_{\B} }$ is center-symmetric.
Provided such an observation, it only remains to show that the associators 
$
(x,y,c)_{\alpha_{\A} \oplus \alpha_{\B} }$,
$(x,b, z)_{\alpha_{\A} \oplus \alpha_{\B} }$,
$(x,b, c)_{\alpha_{\A} \oplus \alpha_{\B} }$, and 
$(a, y,c)_{\alpha_{\A} \oplus \alpha_{\B} }
 $
are center-symmetric.

From the relation \eqref{eq_hom_csa}, we have:
\beqs
(x,y,c)_{\alpha_{\A} \oplus \alpha_{\B} }&=& r_{\B}(\alpha_{\B}(c))(xy)-(\alpha_{\A}(x))\cdot (r_{\B}(c)y)-
r_{\B}(l_{\A}(y)c)(\alpha_{\A}(x))\cr
&+&l_{\A}(x  y)(\alpha_{\B}(c)) -l_{\A}(\alpha_{\A}(x))(l_{\A}(y)c),
\eeqs
\beqs
(x,b, z)_{\alpha_{\A} \oplus \alpha_{\B} }&=&(r_{\B}(b)x)\cdot (\alpha_{\A}(z))+
l_{\B}(l_{\A}(x)b)(\alpha_{\A}(z))-(\alpha_{\A}(x))\cdot (l_{\B}(z))\cr 
&-& r_{\B}(r_{\A}(z)b)(\alpha_{\A}(x))
+ r_{\A}(\alpha_{\A}(z))(l_{\A}(x)b)-l_{\A}(\alpha_{\A}(x))(r_{\A}(z)b),
\eeqs
\beqs
(x,b, c)_{\alpha_{\A} \oplus \alpha_{\B} }&=& r_{\B}(\alpha_{\B}(c))(r_{\B}(b)x)-r_{\B}(b\circ c)(\alpha_{\A}(x))
+(l_{\A}(x)b)\circ (\alpha_{\B}(c))\cr 
&+& l_{\A}(r_{\B}(b)x)(\alpha_{\B}(c))
-l_{\A}(\alpha_{\A}(x))(bc),
\eeqs
\beqs
(a, y,c)_{\alpha_{\A} \oplus \alpha_{\B} }&=& r_{\B}(\alpha_{\B}(c))(l_{\B}(a)y)-
l_{\B}(\alpha_{\B}(a))(r_{\B}(c)y) +(r_{\A}(y)a)(\alpha_{\B}(c))\cr 
&+& l_{\A}(l_{\B}(a)y)(\alpha_{\B}(c))
-\alpha_{\B}(a)(l_{\A}(y)c)-r_{\A}(r_{\B}(c)y)(\alpha_{\B}(a)).
\eeqs
Besides, 
\beqs
(c,y,x)_{\alpha_{\A} \oplus \alpha_{\B} }&=&(l_{\B}(c)y)\cdot (\alpha_{\A}(x))+l_{\B}(r_{\A}(y)c)(\alpha_{\A}(x))
-l_{\B}(\alpha_{\B}(c))(xy)\cr 
&+& r_{\A}(\alpha(x))(r_{\A}(y)c) 
-r_{\A}(yx)(\alpha_{\B}(c)),
\eeqs
\beqs
(z,b, x)_{\alpha_{\A} \oplus \alpha_{\B} }&=& (r_{\B}(b)z)\cdot (\alpha_{\A}(x))+
l_{\B}(l_{\A}(z)b)(\alpha_{\A}(x))-(\alpha_{\A}(z))\cdot (l_{\B}(x))\cr 
&-& r_{\B}(r_{\A}(x)b)(\alpha_{\A}(z))
+r_{\A}(\alpha_{\A}(x))(l_{\A}(z)b)-l_{\A}(\alpha_{\A}(z))(r_{\A}(x)b),
\eeqs
\beqs
(c,b,x)_{\alpha_{\A} \oplus \alpha_{\B} }&=& l_{\B}(c\circ b)\cdot (\alpha_{\A}(x))-l_{\B}(\alpha_{\B}(c))(l_{\B}(b)x)+
r_{\A}(\alpha_{\A}(x))(c\circ b)\cr 
&-&\alpha_{\B}(a)\circ(r_{\A}(x)b) 
-r_{\A}(l_{\B}(b)x)(\alpha_{\B}(c)) ,
\eeqs
\beqs
(c, y,a)_{\alpha_{\A} \oplus \alpha_{\B} }&=& r_{\B}(\alpha_{\B}(a))(l_{\B}(c)y)-
l_{\B}(\alpha_{\B}(c))(r_{\B}(a)y) +(r_{\A}(y)c)(\alpha_{\B}(a))\cr 
&+&l_{\A}(l_{\B}(c)y)(\alpha_{\B}(a))
-\alpha_{\B}(c)(l_{\A}(y)a)-r_{\A}(r_{\B}(a)y)(\alpha_{\B}(c)).
\eeqs
Hence,  $(x,y,c)_{\alpha_{\A} \oplus \alpha_{\B} }=(c,y,x)_{\alpha_{\A} \oplus \alpha_{\B} }$ 
is equivalent to the relations \eqref{eq_matched_hom_center_1} and \eqref{bimod_hom_center2},\\
$(x,b, z)_{\alpha_{\A} \oplus \alpha_{\B} }=(z,b, x)_{\alpha_{\A} \oplus \alpha_{\B} }$  to 
the relations \eqref{eq_matched_hom_center_2} and \eqref{bimod_hom_center3}, 
$(x,b, c)_{\alpha_{\A} \oplus \alpha_{\B} }=(c,b, x)_{\alpha_{\A} \oplus \alpha_{\B} }$  to 
the relations \eqref{eq_matched_hom_center_3} and \eqref{bimod_hom_center2},
and finally,
$(a, y,c)_{\alpha_{\A} \oplus \alpha_{\B} }=(c, y,a)_{\alpha_{\A} \oplus \alpha_{\B} }$  to 
the relations \eqref{eq_matched_hom_center_4} and \eqref{bimod_hom_center3}.

Conversely, if $(\A, \cdot, \alpha_{\A})$ and $(\B, \circ, \alpha_{\B})$ are two hom-center-symmetric subalgebras
of the hom-center-symmetric algebra $(\A\oplus\B, \ast, \alpha_{\A}+\alpha_{\B})$, then,  by a direct computation, 
the linear maps $l_{\A}, r_{\A}: \A \rightarrow \mathfrak{gl}(\B)$ and 
$l_{\B}, r_{\B}: \B\rightarrow \mathfrak{gl}(\A), $ given by
\beqs
x\ast a=l_{\A}(x)a+r_{\B}(a)x, \qquad a\ast x= l_{\B}(a)x+r_{\A}(x)a, \qquad \forall x \in \A, a\in \B,
\eeqs
satisfy the relations \eqref{eq_matched_hom_center_1}, \eqref{eq_matched_hom_center_2}, 
\eqref{eq_matched_hom_center_3} and \eqref{eq_matched_hom_center_4}. In addition, 
$(l_{\A}, r_{\A}, \B)$ and $(l_{\B}, r_{\B}, \A)$ are  bimodules of $(\A, \alpha_{\A})$ and 
$(\B, \alpha_{\B}),$ respectively. 
$\cqfd$

We denote the hom-center-symmetric algebra $(\A\oplus\B, \ast, \alpha_{{\A}}\oplus \alpha_{\B})$ by 
$(\A \bowtie_{l_{\B}, r_{\B}}^{l_{\A}, r_{\A}}\B, \alpha_{\A}\oplus\alpha_{\B}).$ The octuple
$(l_{\A}, r_{\A}, l_{\B},r_{\B}, \alpha_{\A}, \alpha_{\B}, \A, \B),$  where 
$l_{\A}, r_{\A}, l_{\B},r_{\B}, \alpha_{\A}$ and
$\alpha_{\B}$ satisfy the  conditions \eqref{eq_matched_hom_center_1},
\eqref{eq_matched_hom_center_2}, \eqref{eq_matched_hom_center_3}, and \eqref{eq_matched_hom_center_4},  is called 
a matched pair of the hom-center-symmetric algebras $\A$ and $\B$.
\begin{definition}\label{dfn_matched_pair_hom_Lie}
A matched pair of hom-Lie algebras
 $(\G, \h, \rho_{_{\G}}, \rho_{_{\h}}, \varphi_{_{\G}}, \varphi_{_{\h}})$ consists of two
hom-Lie algebras $(\G, [,]_{_{\G}}, \varphi_{_{\G}})$ and $(\h, [,]_{_{\h}}, \varphi_{_{\h}})$,  together with the associated
hom-Lie algebra representations $\rho_{_{\G}}:\G\rightarrow\mathfrak{gl}(\h)$ and
$\rho_{_{\h}}:\h\rightarrow\mathfrak{gl}(\G)$ defined with respect to $\varphi_{_{\h}}$ and $\varphi_{_{\G}},$ respectively, 
satisfying the following relations
\beq\label{eq_Hom_Lie_1}
\rho_{_{\h}}(\varphi_{_{\h}}(a))[x,y]_{_{\G}}&=&[\rho_{_{\h}}(a)(x), \varphi_{_{\G}}(y)]_{_{\G}}+
[\varphi_{_{\G}}(x),\rho_{_{\h}}(a)(y)]_{_{\G}}\cr 
&+&
\rho_{_{\h}}(\rho_{_{\G}}(y)(a))(\varphi_{_{\G}}(x)) 
-\rho_{_{\h}}(\rho_{_{\G}}(x)(a))(\varphi_{_{\G}}(y)),
\eeq
\beq\label{eq_Hom_Lie_2}
\rho_{_{\G}}(\varphi_{_{\G}}(x))[a,b]_{_{\h}}&=&[\rho_{_{\G}}(x)(a), \varphi_{_{\h}}(b)]_{_{\h}}+
[\varphi_{_{\h}}(a),\rho_{_{\G}}(x)(b)]_{_{\h}}\cr 
&+&
\rho_{_{\G}}(\rho_{_{\h}}(b)(x))(\varphi_{_{\h}}(a))
-\rho_{_{\G}}(\rho_{_{\h}}(a)(x))(\varphi_{_{\h}}(b)).
\eeq
\end{definition}
\begin{corollary}
 Let $(\A, \B, l_{\A}, r_{\A}, l_{\B}, r_{\B}, \alpha_{\A}, \alpha_{\B})$ be a matched pair of
 hom-center - symmetric algebras $(\A, \cdot, \alpha_{\A})$ and $(\B, \circ, \alpha_{\B})$. Then,
  $(\G(\A), \G(\B), l_{\A}-r_{\A}, l_{\B}-r_{\B},\alpha_{\A}, \alpha_{\B})$ is a matched pair of hom-Lie algebras $\G(\A)_{\alpha_{\A}}$ and $\G(\B)_{\alpha_{\B}}$.
\end{corollary}
\textbf{Proof}.

From Lemma~\ref{lemma_hom_1}, it comes that $(\rho_{_{\G(\A)}}:=l_{\A}-r_{\A}, \alpha_{\A}, \B)$ and 
 $(\rho_{_{\G(\B)}}:=l_{\B}-r_{\B}, \alpha_{\B}, \A)$ are  linear representations of the 
sub-adjacent Lie algebras $(\G(\A), \alpha_{\A})$ and $(\G(\B), \alpha_{\B}),$ respectively. It  only remains to show that
the linear maps $(\rho_{_{\G}}=\rho_{_{\G(\A)}}, \varphi_{\G}=\alpha_{\A})$ and 
 $(\rho_{_{\h}}=\rho_{_{\G(\B)}}, \varphi_{\h}=\alpha_{\B})$ satisfy the relations \eqref{eq_Hom_Lie_1}
 and
 \eqref{eq_Hom_Lie_2}.
In fact, for all $x,y,z\in \A$ and  $a, b, c\in \B,$ we have:
\beqs
&&[\alpha_{\A}(x)+\alpha_{\B}(a), [y+b, z+c]]+
[\alpha_{\A}(y)+\alpha_{\B}(b), [z+c, x+a]]\cr 
&+&
[\alpha_{\A}(z)+\alpha_{\B}(c), [x+a, y+b]]
=
[\alpha_{\A}(x)+\alpha_{\B}(a), ([y,z]+(l_{\B}-r_{\B})(b)z-(l_{\B}-r_{\B})(c)y)\cr 
&+&
([b,c]+(l_{\A}-r_{\A})(y)c-
(l_{\A}-r_{\A})(z)b)]
+[\alpha_{\A}(y)+\alpha_{\B}(b), ([z,x]+(l_{\B}-r_{\B})(c)x\cr 
&-&(l_{\B}-r_{\B})(a)z)
+([c,a]+(l_{\A}-r_{\A})(z)a-
(l_{\A}-r_{\A})(x)c)]
+[\alpha_{\A}(z)+\alpha_{\B}(c),
 ([x,y]
 \cr 
 &+&(l_{\B}-r_{\B})(a)y-(l_{\B}-r_{\B})(b)x)+([a,b]+(l_{\A}-r_{\A})(x)b-
(l_{\A}-r_{\A})(y)a)]=
[\alpha_{\A}(x),[y,z]\cr 
&+&(l_{\B}-r_{\B})(b)z-(l_{\B}-r_{\B})(c)y]
+(l_{\B}-r_{\B})(\alpha_{\B}(a))( [y,z]+(l_{\B}-r_{\B})(b)z-(l_{\B}-r_{\B})(c)y)\cr 
&-&(l_{\B}-r_{\B})([b,c]
+(l_{\A}-r_{\A})(y)c-(l_{\A}-r_{\A})(z)b)(\alpha_{\A}(x))+
[\alpha_{\B}(a),[b,c]+(l_{\A}-r_{\A})(y)c \cr 
&-&(l_{\A}-r_{\A})(z)b]
+(l_{\A}-r_{\A})(\alpha_{\A}(x))( [b,c]+(l_{\A}-r_{\A})(y)c-(l_{\A}-r_{\A})(z)b)-(l_{\A}-r_{\A})([y,z]
\cr 
&+&(l_{\B}-r_{\B})(b)z
-(l_{\B}-r_{\B})(c)y)(\alpha_{\B}(a))+[\alpha_{\A}(y),[z,x]+(l_{\B}-r_{\B})(c)x
-(l_{\B}-r_{\B})(a)z]\cr 
&+&(l_{\B}-r_{\B})(\alpha_{\B}(b))( [z,x]+(l_{\B}-r_{\B})(c)x-(l_{\B}-r_{\B})(a)z)
-(l_{\B}-r_{\B})([c,a]+(l_{\A}-r_{\A})(z)a\cr 
&-&(l_{\A}-r_{\A})(x)c)(\alpha_{\A}(y))
+[\alpha_{\B}(b),[c,a]
+(l_{\A}-r_{\A})(z)a-(l_{\A}-r_{\A})(x)c]\cr 
&+&(l_{\A}-r_{\A})(\alpha_{\A}(y))( [c,a]+(l_{\A}-r_{\A})(z)a
-(l_{\A}-r_{\A})(x)c)-(l_{\A}-r_{\A})([z,x]+(l_{\B}-r_{\B})(c)x\cr 
&-&(l_{\B}-r_{\B})(a)z)(\alpha_{\B}(b))
+[\alpha_{\A}(z),[x,y]+(l_{\B}-r_{\B})(a)y-(l_{\B}-r_{\B})(b)x]
\cr 
&+&(l_{\B}-r_{\B})(\alpha_{\B}(c))( [x,y]
+(l_{\B}-r_{\B})(a)y-(l_{\B}-r_{\B})(b)x)
-(l_{\B}-r_{\B})([a,b]+(l_{\A}-r_{\A})(x)b
\cr 
&-&(l_{\A}-r_{\A})(y)a)(\alpha_{\A}(z))+[\alpha_{\B}(c),[a,b]+(l_{\A}-r_{\A})(x)b-(l_{\A}-r_{\A})(y)a]\cr 
&+&(l_{\A}-r_{\A})(\alpha_{\A}(z))( [a,b]+(l_{\A}-r_{\A})(x)b-(l_{\A}-r_{\A})(y)a)-(l_{\A}-r_{\A})([x,y]
+(l_{\B}-r_{\B})(a)y\cr 
&-&(l_{\B}-r_{\B})(b)x)(\alpha_{\B}(c))=[\alpha_{\A}(x),[y,z]]+
[\alpha_{\A}(x), (l_{\B}-r_{\B})(b)z]-
[\alpha_{\A}(x),(l_{\B}-r_{\B})(c)y]\cr 
&+&(l_{\B}-r_{\B})(\alpha_{\B}(a))( [y,z])+
(l_{\B}-r_{\B})(\alpha_{\B}(a))((l_{\B}-r_{\B})(b)z)
-(l_{\B}-r_{\B})(\alpha_{\B}(a))((l_{\B}-r_{\B})(c)y\cr 
&-&(l_{\B}-r_{\B})([b,c])(\alpha_{\A}(x))-(l_{\B}-r_{\B})((l_{\A}-r_{\A})(y)c)(\alpha_{\A}(x))\cr 
&+&(l_{\B}-r_{\B})((l_{\A}-r_{\A})(z)b)(\alpha_{\A}(x))+[\alpha_{\B}(a),[b,c]]+
[\alpha_{\B}(a),(l_{\A}-r_{\A})(y)c]\cr 
&-&
[\alpha_{\B}(a),(l_{\A}-r_{\A})(z)b]+(l_{\A}-r_{\A})(\alpha_{\A}(x))( [b,c])+
(l_{\A}-r_{\A})(\alpha_{\A}(x))((l_{\A}-r_{\A})(y)c)\cr 
&-&
(l_{\A}-r_{\A})(\alpha_{\A}(x))((l_{\A}-r_{\A})(z)b)-(l_{\A}-r_{\A})([y,z])(\alpha_{\B}(a))\cr 
&-&
(l_{\A}-r_{\A})( (l_{\B}-r_{\B})(b)z)(\alpha_{\B}(a))+(l_{\A}-r_{\A})((l_{\B}-r_{\B})(c)y)(\alpha_{\B}(a))+[\alpha_{\A}(y),[z,x]]
\cr&+& 
[\alpha_{\A}(y),(l_{\B}-r_{\B})(c)x]- 
[\alpha_{\A}(y),(l_{\B}-r_{\B})(a)z]+(l_{\B}-r_{\B})(\alpha_{\B}(b))( [z,x])\cr 
&+&
(l_{\B}-r_{\B})(\alpha_{\B}(b))((l_{\B}-r_{\B})(c)x)
-(l_{\B}-r_{\B})(\alpha_{\B}(b))( (l_{\B}-r_{\B})(a)z)\cr 
&-&(l_{\B}-r_{\B})([c,a])(\alpha_{\A}(y)) 
-(l_{\B}-r_{\B})((l_{\A}-r_{\A})(z)a)(\alpha_{\A}(y))\cr 
&+&(l_{\B}-r_{\B})((l_{\A}-r_{\A})(x)c)(\alpha_{\A}(y))+[\alpha_{\B}(b),[c,a]]+
[\alpha_{\B}(b),(l_{\A}-r_{\A})(z)a]\cr 
&-&
[\alpha_{\B}(b),(l_{\A}-r_{\A})(x)c]+(l_{\A}-r_{\A})(\alpha_{\A}(y))([c,a])+ 
(l_{\A}-r_{\A})(\alpha_{\A}(y))((l_{\A}-r_{\A})(z)a) \cr 
&-& 
(l_{\A}-r_{\A})(\alpha_{\A}(y))((l_{\A}-r_{\A})(x)c)-(l_{\A}-r_{\A})([z,x])(\alpha_{\B}(b))\cr
&-&
(l_{\A}-r_{\A})((l_{\B}-r_{\B})(c)x)(\alpha_{\B}(b))+(l_{\A}-r_{\A})((l_{\B}-r_{\B})(a)z)(\alpha_{\B}(b))
+[\alpha_{\A}(z),[x,y]]\cr 
&+&
[\alpha_{\A}(z),(l_{\B}-r_{\B})(a)y]-
[\alpha_{\A}(z),(l_{\B}-r_{\B})(b)x]+(l_{\B}-r_{\B})(\alpha_{\B}(c))([x,y])\cr 
&+&
(l_{\B}-r_{\B})(\alpha_{\B}(c))((l_{\B}-r_{\B})(a)y)
-(l_{\B}-r_{\B})(\alpha_{\B}(c))((l_{\B}-r_{\B})(b)x)\cr 
&-&(l_{\B}-r_{\B})([a,b])(\alpha_{\A}(z))
-(l_{\B}-r_{\B})((l_{\A}-r_{\A})(x)b)(\alpha_{\A}(z))\cr 
&+&(l_{\B}-r_{\B})((l_{\A}-r_{\A})(y)a)(\alpha_{\A}(z))+[\alpha_{\B}(c),[a,b]]+ 
[\alpha_{\B}(c),(l_{\A}-r_{\A})(x)b]\cr 
&-&
[\alpha_{\B}(c),(l_{\A}-r_{\A})(y)a]+(l_{\A}-r_{\A})(\alpha_{\A}(z))( [a,b])+ 
(l_{\A}-r_{\A})(\alpha_{\A}(z))((l_{\A}-r_{\A})(x)b)\cr
&-&
(l_{\A}-r_{\A})(\alpha_{\A}(z))((l_{\A}-r_{\A})(y)a)-(l_{\A}-r_{\A})([x,y])(\alpha_{\B}(c))\cr 
&-&
(l_{\A}-r_{\A})((l_{\B}-r_{\B})(a)y)(\alpha_{\B}(c))+
(l_{\A}-r_{\A})((l_{\B}-r_{\B})(b)x)(\alpha_{\B}(c))=[\alpha_{\A}(x),[y,z]]\cr 
&+&[\alpha_{\A}(y),[z,x]]+[\alpha_{\A}(z),[x,y]] 
+[\alpha_{\B}(a),[b,c]]+[\alpha_{\B}(b),[c,a]]+[\alpha_{\B}(c),[a,b]]\cr 
&+&[\alpha_{\A}(x), (l_{\B}-r_{\B})(b)z]
+(l_{\B}-r_{\B})((l_{\A}-r_{\A})(z)b)(\alpha_{\A}(x))+(l_{\B}-r_{\B})(\alpha_{\B}(b))( [z,x])\cr 
&-&[\alpha_{\A}(z),(l_{\B}-r_{\B})(b)x]- 
(l_{\B}-r_{\B})((l_{\A}-r_{\A})(x)b)(\alpha_{\A}(z))+[\alpha_{\A}(y),(l_{\B}-r_{\B})(c)x]\cr 
&-&[\alpha_{\A}(x),(l_{\B}-r_{\B})(c)y]+(l_{\B}-r_{\B})(\alpha_{\B}(c))([x,y])
-(l_{\B}-r_{\B})((l_{\A}-r_{\A})(y)c)(\alpha_{\A}(x))\cr 
&+&(l_{\B}-r_{\B})((l_{\A}-r_{\A})(x)c)(\alpha_{\A}(y))+(l_{\B}-r_{\B})(\alpha_{\B}(a))( [y,z])-
 [\alpha_{\A}(y),(l_{\B}-r_{\B})(a)z]\cr
&+& 
(l_{\B}-r_{\B})((l_{\A}-r_{\A})(y)a)(\alpha_{\A}(z))
-(l_{\B}-r_{\B})((l_{\A}-r_{\A})(z)a)(\alpha_{\A}(y))\cr 
&+&[\alpha_{\A}(z),(l_{\B}-r_{\B})(a)y]
+(l_{\B}-r_{\B})(\alpha_{\B}(a))((l_{\B}-r_{\B})(b)z)-(l_{\B}-r_{\B})(\alpha_{\B}(b))( (l_{\B}-r_{\B})(a)z)\cr
&-&
(l_{\B}-r_{\B})([a,b])(\alpha_{\A}(z))
+(l_{\B}-r_{\B})(\alpha_{\B}(c))((l_{\B}-r_{\B})(a)y)\cr 
&
-&(l_{\B}-r_{\B})(\alpha_{\B}(a))((l_{\B}-r_{\B})(c)y)
-
(l_{\B}-r_{\B})([c,a])(\alpha_{\A}(y))\cr 
&+&(l_{\B}-r_{\B})(\alpha_{\B}(b))((l_{\B}-r_{\B})(c)x)-
[\alpha_{\B}(c),(l_{\A}-r_{\A})(y)a]\cr 
&-&(l_{\B}-r_{\B})([b,c])(\alpha_{\A}(x))-
(l_{\B}-r_{\B})(\alpha_{\B}(c))((l_{\B}-r_{\B})(b)x)
+[\alpha_{\B}(a),(l_{\A}-r_{\A})(y)c]\cr 
&+& (l_{\A}-r_{\A})((l_{\B}-r_{\B})(c)y)(\alpha_{\B}(a))
+(l_{\A}-r_{\A})(\alpha_{\A}(y))([c,a])\cr  
&-&(l_{\A}-r_{\A})((l_{\B}-r_{\B})(a)y)(\alpha_{\B}(c))+[\alpha_{\B}(b),(l_{\A}-r_{\A})(z)a]
-[\alpha_{\B}(a),(l_{\A}-r_{\A})(z)b]\cr 
&-&(l_{\A}-r_{\A})( (l_{\B}-r_{\B})(b)z)(\alpha_{\B}(a))
+(l_{\A}-r_{\A})((l_{\B}-r_{\B})(a)z)(\alpha_{\B}(b))\cr 
&+&(l_{\A}-r_{\A})(\alpha_{\A}(z))( [a,b])
+(l_{\A}-r_{\A})(\alpha_{\A}(x))([b,c])-[\alpha_{\B}(b),(l_{\A}-r_{\A})(x)c]\cr 
&-&
(l_{\A}-r_{\A})((l_{\B}-r_{\B})(c)x)(\alpha_{\B}(b))
+(l_{\A}-r_{\A})((l_{\B}-r_{\B})(b)x)(\alpha_{\B}(c))+[\alpha_{\B}(c),(l_{\A}-r_{\A})(x)b]\cr 
&+&(l_{\A}-r_{\A})(\alpha_{\A}(x))((l_{\A}-r_{\A})(y)c)- (l_{\A}-r_{\A})(\alpha_{\A}(y))((l_{\A}-r_{\A})(x)c\cr 
&-&(l_{\A}-r_{\A})([x,y])(\alpha_{\B}(c))+(l_{\A}-r_{\A})(\alpha_{\A}(y))((l_{\A}-r_{\A})(z)a)
-(l_{\A}-r_{\A})([y,z])(\alpha_{\B}(a))\cr 
&-&
(l_{\A}-r_{\A})(\alpha_{\A}(z))((l_{\A}-r_{\A})(y)a)
+(l_{\A}-r_{\A})(\alpha_{\A}(z))((l_{\A}-r_{\A})(x)b)
\cr
&-&
(l_{\A}-r_{\A})([z,x])(\alpha_{\B}(b))-(l_{\A}-r_{\A})(\alpha_{\A}(x))((l_{\A}-r_{\A})(z)b)= 
[\alpha_{\A}(x), (l_{\B}-r_{\B})(b)z]\cr 
&+&(l_{\B}-r_{\B})((l_{\A}-r_{\A})(z)b)(\alpha_{\A}(x))+(l_{\B}-r_{\B})(\alpha_{\B}(b))( [z,x])-[\alpha_{\A}(z),(l_{\B}-r_{\B})(b)x]\cr 
&-& 
(l_{\B}-r_{\B})((l_{\A}-r_{\A})(x)b)(\alpha_{\A}(z))+[\alpha_{\A}(y),(l_{\B}-r_{\B})(c)x]\cr 
&+&(l_{\B}-r_{\B})((l_{\A}-r_{\A})(x)c)(\alpha_{\A}(y))+
(l_{\B}-r_{\B})(\alpha_{\B}(c))([x,y])
\cr 
&-&(l_{\B}-r_{\B})((l_{\A}-r_{\A})(y)c)(\alpha_{\A}(x))-[\alpha_{\A}(x),(l_{\B}-r_{\B})(c)y]+(l_{\B}-r_{\B})(\alpha_{\B}(a))( [y,z]) \cr 
&-& [\alpha_{\A}(y),(l_{\B}-r_{\B})(a)z]+ 
(l_{\B}-r_{\B})((l_{\A}-r_{\A})(y)a)(\alpha_{\A}(z))\cr  
&-&(l_{\B}-r_{\B})((l_{\A}-r_{\A})(z)a)(\alpha_{\A}(y))+[\alpha_{\A}(z),(l_{\B}-r_{\B})(a)y]
+[\alpha_{\B}(a),(l_{\A}-r_{\A})(y)c]\cr 
&+&(l_{\A}-r_{\A})((l_{\B}-r_{\B})(c)y)(\alpha_{\B}(a))-[\alpha_{\B}(c),(l_{\A}-r_{\A})(y)a]
+(l_{\A}-r_{\A})(\alpha_{\A}(y))([c,a])\cr 
&-& (l_{\A}-r_{\A})((l_{\B}-r_{\B})(a)y)(\alpha_{\B}(c))+[\alpha_{\B}(b),(l_{\A}-r_{\A})(z)a]
-[\alpha_{\B}(a),(l_{\A}-r_{\A})(z)b]\cr 
&-&(l_{\A}-r_{\A})( (l_{\B}-r_{\B})(b)z)(\alpha_{\B}(a))+(l_{\A}-r_{\A})((l_{\B}-r_{\B})(a)z)(\alpha_{\B}(b))\cr 
&+&(l_{\A}-r_{\A})(\alpha_{\A}(z))( [a,b])
+(l_{\A}-r_{\A})(\alpha_{\A}(x))([b,c])-[\alpha_{\B}(b),(l_{\A}-r_{\A})(x)c]\cr 
&-&
(l_{\A}-r_{\A})((l_{\B}-r_{\B})(c)x)(\alpha_{\B}(b))+(l_{\A}-r_{\A})((l_{\B}-r_{\B})(b)x)(\alpha_{\B}(c))+
[\alpha_{\B}(c),(l_{\A}-r_{\A})(x)b].
\eeqs
Therefore, the hom-Jacobi identity \eqref{eq_hom_Jacobi_identity} condition associated to the underlying
hom-Lie algebra from the quadruple  $(\G(\A), \G(\B), l_{\A}-r_{\A}, l_{\B}-r_{\B},\alpha_{\A}, \alpha_{\B})$
is equivalent to the relations \eqref{eq_Hom_Lie_1} and \eqref{eq_Hom_Lie_2}.
$\cqfd$
\begin{theorem}\label{Theo_Hom_ll}
Let $(\A, \cdot, \alpha)$ be a hom-center-symmetric algebra. Suppose there is a 
hom-center-symmetric algebra structure $"\circ "$ on its dual vector space $\A^*$. 
  $(\A, \A^*, R_{\cdot}^*, L_{\cdot}^*, R_{\circ}^*, L_{\circ}^*, \alpha, \alpha^*)$ 
is a matched pair of hom-center-symmetric algebras $(\A, \cdot, \alpha)$ and $(\A^*, \circ, \alpha^*)$ if and only
if the sixtuple $(\G(\A),\G(\A^*),-\ad_{\cdot}^*=R_{\cdot}^*-L_{\cdot}^*,-\ad_{\circ}^*=R_{\circ}^*-L_{\circ}^*,\alpha,\alpha^*)$
 is a matched pair
of the hom-Lie algebras $\G(\A)_{\alpha}$ and $\G(\A^*)_{\alpha^*}$.
\end{theorem}
\textbf{Proof}.

By considering  Theorem~\ref{thm_matched_hom_center},  setting $\B=\A^*$, $l_{\A}=R_{\cdot}^*$, $r_{\A}=L_{\cdot}^*$,
$l_{\B}=R_{\circ}^*$, $r_{\B}=L_{\circ}^*$, using  Definition~\ref{dfn_matched_pair_hom_Lie} by assuming
that $\G(\A)=\G$, $\G(\A^*)=\h$, $\rho_{_{\G}}=R_{\cdot}^*-L_{\cdot}^*$, $\rho_{_{\h}}=R_{\circ}^*-L_{\circ}^*$, 
$\varphi_{_{\G}}=\alpha$, $\varphi_{_{\h}}=\alpha^*$, and taking into account the relations \eqref{dual_linear_1}
 and \eqref{dual_linear_2},
we get the following equivalences:
\begin{itemize}
\item The equation \eqref{eq_Hom_Lie_1} is equivalent to both the equations \eqref{eq_matched_hom_center_1} and 
\eqref{eq_matched_hom_center_2},  i.e. for all $x,y\in \A$ and  $a\in \A^*,$ we have:
\beqs
&&-(R_{\circ}^*-L_{\circ}^*)(\alpha^*(a))[x,y]_{_{\G}}+[(R_{\circ}^*-L_{\circ}^*)(a)(x), \alpha(y)]_{_{\G}}+
[\alpha(x),(R_{\circ}^*-L_{\circ}^*)(a)(y)]_{_{\G}}\cr 
&&+(R_{\circ}^*-L_{\circ}^*)( (R_{\cdot}^*-L_{\cdot}^*) (y)(a))(\alpha(x))- 
(R_{\circ}^*-L_{\circ}^*)( (R_{\cdot}^*-L_{\cdot}^*)(x)(a))(\alpha(y))= \cr 
&&-R_{\circ}^* (\alpha^*(a))(xy) +R_{\circ}^* (\alpha^*(a))(yx)
+L_{\circ}^*(\alpha^*(a))(xy)-L_{\circ}^*(\alpha^*(a))(yx) 
+(R_{\circ}^*(a)(x))\cdot (\alpha(y))\cr 
&&-(L_{\circ}^*(a)(x))\cdot (\alpha(y))
-(\alpha(y))\cdot(R_{\circ}^*(a)(x))+ (\alpha(y))\cdot(L_{\circ}^*(a)(x)) 
-(R_{\circ}^*(a)(y))\cdot (\alpha(x))\cr 
&&+(L_{\circ}^*(a)(y))\cdot (\alpha(x))+(\alpha(x))\cdot(R_{\circ}^*(a)(y))- (\alpha(x))\cdot(L_{\circ}^*(a)(y))
+R_{\circ}^*( (R_{\cdot}^*) (y)(a))(\alpha(x))\cr 
&&- R_{\circ}^*((L_{\cdot}^*) (y)(a)) (\alpha(x))
-L_{\circ}^*( (R_{\cdot}^* ) (y)(a)) (\alpha(x))
+L_{\circ}^*( ( L_{\cdot}^*) (y)(a)) (\alpha(x)) \cr 
&& -R_{\circ}^*( (R_{\cdot}^*) (x)(a))(\alpha(y))+R_{\circ}^*((L_{\cdot}^*) (x)(a)) (\alpha(y))
+L_{\circ}^*( (R_{\cdot}^* ) (x)(a)) (\alpha(y))\cr 
&&
-L_{\circ}^*( ( L_{\cdot}^*) (x)(a)) (\alpha(y))
=\{-R_{\circ}^* (\alpha^*(a))(xy)-L_{\circ}^*(\alpha^*(a))(yx)+(R_{\circ}^*(a)(x))\cdot (\alpha(y))\cr 
&&+
(\alpha(y))\cdot(L_{\circ}^*(a)(x))+ L_{\circ}^*( (R_{\cdot}^* ) (x)(a)) (\alpha(y))+
R_{\circ}^*((L_{\cdot}^*) (x)(a)) (\alpha(y)) \}\cr 
&&+
\{R_{\circ}^* (\alpha^*(a))(yx)+L_{\circ}^*(\alpha^*(a))(xy)- (R_{\circ}^*(a)(y))\cdot (\alpha(x))-L_{\circ}^*( (R_{\cdot}^* ) (y)(a)) (\alpha(x))\cr 
&&- 
(\alpha(x))\cdot(L_{\circ}^*(a)(y)) - R_{\circ}^*((L_{\cdot}^*) (y)(a)) (\alpha(x))\}+\{ 
-(L_{\circ}^*(a)(x))\cdot (\alpha(y))
-(\alpha(y))\cdot(R_{\circ}^*(a)(x))\cr 
&& 
-  R_{\circ}^*( (R_{\cdot}^*) (x)(a))(\alpha(y))
-L_{\circ}^*( ( L_{\cdot}^*) (x)(a)) (\alpha(y))
+(L_{\circ}^*(a)(y))\cdot (\alpha(x))+(\alpha(x))\cdot(R_{\circ}^*(a)(y))\cr 
&&
+R_{\circ}^*( (R_{\cdot}^*) (y)(a))(\alpha(x))
+L_{\circ}^*(  L_{\cdot}^*(y)(a)) (\alpha(x))
\}=0,
\eeqs
using the fact that the two first relations give zero from the equation \eqref{eq_matched_hom_center_1}, 
and the last  brace  yields zero from the equation \eqref{eq_matched_hom_center_2}.
\item The equation \eqref{eq_Hom_Lie_2} is equivalent to both the equations \eqref{eq_matched_hom_center_3} and 
\eqref{eq_matched_hom_center_4}, i.e. for all $a, b\in \A^*$ and  $x\in \A$, we have:
\beqs
&&-(R_{\cdot}^*-L_{\cdot}^*)(\alpha(x))[a,b]_{_{\h}}+[(R_{\cdot}^*-L_{\cdot}^*)(x)(a), \alpha^*(b)]_{_{\h}}+
[\alpha^*(a),(R_{\cdot}^*-L_{\cdot}^*)(x)(b)]_{_{\h}}\cr 
&&+(R_{\cdot}^*-L_{\cdot}^*)( (R_{\circ}^*-L_{\circ}^*) (b)(x))(\alpha^*(a))- 
(R_{\cdot}^*-L_{\cdot}^*)( (R_{\circ}^*-L_{\circ}^*)(a)(x))(\alpha^*(b))= \cr 
&&-R_{\cdot}^* (\alpha(x))(a\circ b) +R_{\cdot}^* (\alpha(x))(b\circ a)
+L_{\cdot}^*(\alpha(x))(a\circ b)-L_{\cdot}^*(\alpha(x))(b\circ a)\cr 
&& 
+(R_{\cdot}^*(x)(a))\circ (\alpha^*(b)) -(L_{\cdot}^*(x)(a))\circ (\alpha^*(b))
-(\alpha^*(b))\circ(R_{\cdot}^*(x)(a))+ (\alpha^*(b))\circ(L_{\cdot}^*(x)(a))\cr 
&& 
-(R_{\cdot}^*(x)(b))\circ (\alpha^*(a))+(L_{\cdot}^*(x)(b))\circ (\alpha^*(a))+(\alpha^*(a))\circ(R_{\cdot}^*(x)(b))- (\alpha^*(a))\circ(L_{\cdot}^*(x)(b))\cr 
&&
+R_{\cdot}^*( (R_{\circ}^*) (b)(x))(\alpha^*(a))- R_{\cdot}^*((L_{\circ}^*) (b)(x)) (\alpha^*(a))
-L_{\cdot}^*( (R_{\circ}^* ) (b)(x)) (\alpha^*(a))\cr 
&&
+L_{\cdot}^*( ( L_{\circ}^*) (b)(x)) (\alpha^*(a))- R_{\cdot}^*( (R_{\circ}^*) (a)(x))(\alpha^*(b))+R_{\cdot}^*((L_{\circ}^*) (a)(x)) (\alpha^*(b))\cr 
&&
+L_{\cdot}^*( (R_{\circ}^* ) (a)(x)) (\alpha^*(b))
-L_{\cdot}^*( ( L_{\circ}^*) (a)(x)) (\alpha^*(b))
=\{-R_{\cdot}^* (\alpha(x))(a\circ b)-L_{\cdot}^*(\alpha(x))(b\circ a)\cr 
&& +(R_{\cdot}^*(x)(a))\circ (\alpha^*(b)) 
+(\alpha^*(b))\circ(L_{\cdot}^*(x)(a))+ L_{\cdot}^*( (R_{\circ}^* ) (a)(x)) (\alpha^*(b))\cr 
&& +
R_{\cdot}^*((L_{\circ}^*) (a)(x)) (\alpha^*(b)) \}+
\{R_{\cdot}^* (\alpha(x))(b\circ a)+ L_{\cdot}^*(\alpha(x))(a\circ b) -    (R_{\cdot}^*(x)(b))\circ (\alpha^*(a))\cr 
&&-L_{\cdot}^*( R_{\circ}^*  (b)(x)) (\alpha^*(a))- 
(\alpha^*(a))\circ(L_{\cdot}^*(x)(b)) - R_{\cdot}^*((L_{\circ}^*) (b)(x)) (\alpha^*(a))\}\cr 
&&+\{ 
-(L_{\cdot}^*(x)(a))\circ (\alpha^*(b))
-(\alpha^*(b))\circ(R_{\cdot}^*(x)(a))  - R_{\cdot}^*( (R_{\circ}^*) (a)(x))(\alpha^*(b))\cr 
&& 
-L_{\cdot}^*( ( L_{\circ}^*) (a)(x)) (\alpha^*(b))
+(L_{\cdot}^*(x)(b))\circ (\alpha^*(a))+(\alpha^*(a))\circ(R_{\cdot}^*(x)(b))\cr 
&&
+R_{\cdot}^*( (R_{\circ}^*) (b)(x))(\alpha^*(a))
+L_{\cdot}^*(  L_{\circ}^*(b)(x)) (\alpha^*(a))
\}=0
\eeqs
because  the two first relations give zero from the equation \eqref{eq_matched_hom_center_3},  
and the last  brace  yields zero from the equation \eqref{eq_matched_hom_center_4}.
\end{itemize}
Therefore,  the equivalence is obtained. 
$\cqfd$
\section{Manin triple and hom-center-symmetric bialgebras}
This section is devoted to the construction of  hom-center-symmetric bialgebras and to  basic properties of the Manin triple of   hom-center-symmetric algebras. 
\begin{theorem}
Consider a hom-Lie algebra $(\G, [,]_{_{\G}}, \alpha_{_{\G}}),$ and its   two   representations 
$(\rho_{_U}, U, \varphi_{_U})$
and $(\rho_{_V}, V, \varphi_{_V})$ on the vector spaces $U$ and $V,$  respectively. Let 
$  \rho_{_U}: \G\rightarrow \mathfrak{gl}(U)$
$\rho_{_V}: \G\rightarrow \mathfrak{gl}(V)$, $\varphi_{_U}: U\rightarrow U$ and $\varphi_{_V}: V\rightarrow V$
be four linear maps
 such that $(\rho_{_U}, \varphi)$ and $(\rho_{_V}, \varphi_{_V})$ satisfy the relations
\eqref{eq_rep_Hom_Lie1} and \eqref{eq_rep_Hom_Lie2}. Then, the linear map $\rho_{_U}\otimes \varphi_{_V}+\varphi_{_U} \otimes\rho_{_V}:
\G\rightarrow \mathfrak{gl}(U\otimes V),$ defined by:  
$(\rho_{_U}\otimes \varphi_{_V}+\varphi_{_U} \otimes\rho_{_V})(x)(u\otimes v)=\rho_{_U}(x)u \otimes \varphi_{_V}+\varphi_{_U} \otimes\rho_{_V}(x)v$ for all
$x\in \G$ , $u\in U$ and $v\in V,$ is  a representation of the hom-Lie algebra $(\G, [,]_{_{\G}}, \alpha_{_{\G}})$
on  the vector space $U\otimes V$.
\end{theorem}
\textbf{Proof}.

Let $(\G, [,]_{_{\G}}, \alpha_{_{\G}})$ be a hom-Lie algebra, $(\rho_{_U}, U, \varphi_{_U})$ and 
$(\rho_{_V}, V, \varphi_{_V})$ be two representations of $(\G, [,]_{_{\G}}, \alpha_{_{\G}})$ on 
$U$ and $V,$ respectively. For all $x\in\G$, $u\in U$ and $v\in V,$ we have: 
\beqs
&&(\rho_{_U}\otimes \varphi_{_V}+\varphi_{_U} \otimes\rho_{_V})(\alpha_{_{\G}}(x))\circ (\varphi_{_U}\otimes \varphi_{_V})
=
(\rho_{_U}(\alpha_{_{\G}}(x))\circ (\varphi_{_U})\otimes (\varphi_{_V}\circ  \varphi_{_V}) \cr 
&&+(\varphi_{_U}\circ  \varphi_{_U}) \otimes(\rho_{_V}(\alpha_{_{\G}}(x))\circ (\varphi_{_U})) 
= 
(\varphi_{_U}\circ \rho_{_U}(x))\otimes (\varphi_{_V}\circ  \varphi_{_V})) 
+(\varphi_{_U}\circ  \varphi_{_U}))\otimes (\varphi_{_V}\circ \rho_{_V}(x))\cr 
&&=
(\varphi_{_U}\otimes \varphi_{_V})\circ( \rho_{_U}(x)\otimes \varphi_{_V}  ) 
+(\varphi_{_U}\otimes \varphi_{_V})\circ(\varphi_{_U} \otimes  \rho_{_V}(x)  ) 
\eeqs
\beqs
(\rho_{_U}\otimes \varphi_{_V}+\varphi_{_U} \otimes\rho_{_V})(\alpha_{_{\G}}(x))\circ (\varphi_{_U}\otimes \varphi_{_V})
=
(\varphi_{_U}\otimes \varphi_{_V})\circ( \rho_{_U}\otimes \varphi_{_V}+\varphi_{_U} \otimes\rho_{_V} ).
\eeqs
In addition, we have
\beqs
&&(\rho_{_U}\otimes \varphi_{_V}+\varphi_{_U}\otimes \rho_{_V})(\alpha_{_{\G}}(x))(\rho_{_U}\otimes \varphi_{_V}+\varphi_{_U}\otimes \rho_{_V})(y)
\cr 
&
-&(\rho_{_U}\otimes \varphi_{_V}+\varphi_{_U}\otimes \rho_{_V})(\alpha_{_{\G}}(y))(\rho_{_U}\otimes \varphi_{_V}+\varphi_{_U}\otimes \rho_{_V})(x)\cr
&=&(\rho_{_U}(\alpha_{_{\G}}(x))\otimes \varphi_{_V}+\varphi_{_U}\otimes \rho_{_V}(\alpha_{_{\G}}(x)))(\rho_{_U}(y)\otimes \varphi_{_V}+
\varphi_{_U}\otimes \rho_{_V}(y))\cr 
&-&
(\rho_{_U}(\alpha_{_{\G}}(y))\otimes \varphi_{_V}+\varphi_{_U}\otimes \rho_{_V}(\alpha_{_{\G}}(y)))(\rho_{_U}(x)\otimes \varphi_{_V}+
\varphi_{_U}\otimes \rho_{_V}(x))\cr
&=&(\rho_{_U}(\alpha_{_{\G}}(x))\otimes \varphi_{_V})(\rho_{_U}(y)\otimes \varphi_{_V})+
(\rho_{_U}(\alpha_{_{\G}}(x))\otimes \varphi_{_V})(\varphi_{_U}\otimes \rho_{_V}(y))
\cr 
&+&(\varphi_{_U}\otimes \rho_{_V}(\alpha_{_{\G}}(x)))(\rho_{_U}(y)\otimes \varphi_{_V})+
(\varphi_{_U}\otimes \rho_{_V}(\alpha_{_{\G}}(x)))(\varphi_{_U}\otimes \rho_{_V}(y))\cr
&-&(\rho_{_U}(\alpha_{_{\G}}(y))\otimes \varphi_{_V})(\rho_{_U}(x)\otimes \varphi_{_V})-
(\rho_{_U}(\alpha_{_{\G}}(y))\otimes \varphi_{_V})(\varphi_{_U}\otimes \rho_{_V}(x))\cr
&-&(\varphi_{_U}\otimes \rho_{_V}(\alpha_{_{\G}}(y)))(\rho_{_U}(x)\otimes \varphi_{_V})-
(\varphi_{_U}\otimes \rho_{_V}(\alpha_{_{\G}}(y)))(\varphi_{_U}\otimes \rho_{_V}(x))\cr
&=&(\rho_{_U}(\alpha_{_{\G}}(x))(\rho_{_U}(y)))\otimes (\varphi_{_V}\circ  \varphi_{_V})+
(\rho_{_U}(\alpha_{_{\G}}(x))\varphi_{_U})\otimes (\varphi_{_V}\rho_{_V}(y))\cr 
&+&(\varphi_{_U}\rho_{_U}(y)\otimes \rho_{_V}(\alpha_{_{\G}}(x))) \varphi_{_V})+
(\varphi_{_U}\circ \varphi_{_U})\otimes (\rho_{_V}(\alpha_{_{\G}}(x))\rho_{_V}(y))\cr
&-&(\rho_{_U}(\alpha_{_{\G}}(y))\rho_{_U}(x))\otimes (\varphi_{_V}\circ \varphi_{_V})-
(\rho_{_U}(\alpha_{_{\G}}(y))\varphi_{_U})\otimes (\varphi_{_V}\rho_{_V}(x))\cr
&-&(\varphi_{_U}\rho_{_U}(x)\otimes (\rho_{_V}(\alpha_{_{\G}}(y)) \varphi_{_V})-
(\varphi_{_U}\circ \varphi_{_U})\otimes (\rho_{_V}(\alpha_{_{\G}}(y)) \rho_{_V}(x))\cr
&=&(\rho_{_U}(\alpha_{_{\G}}(x))(\rho_{_U}(y)))\otimes (\varphi_{_V}\circ  \varphi_{_V})+
(\varphi_{_U}\rho_{_U}(x))\otimes(\varphi_{_V}\rho_{_V}(y))\cr 
&+& (\varphi_{_U}\rho_{_U}(y))\otimes  (\varphi_{_V}\rho_{_V}(x))+
(\varphi_{_U}\circ \varphi_{_U})\otimes (\rho_{_V}(\alpha_{_{\G}}(x))\rho_{_V}(y))\cr
&-&(\rho_{_U}(\alpha_{_{\G}}(y))\rho_{_U}(x))\otimes (\varphi_{_V}\circ \varphi_{_V})-
(\varphi_{_U}\rho_{_U}(y))\otimes (\varphi_{_V}\rho_{_V}(x))\cr
&-&(\varphi_{_U}\rho_{_U}(x))\otimes (\varphi_{_V}\rho_{_V}(y) )-
(\varphi_{_U}\circ \varphi_{_U})\otimes (\rho_{_V}(\alpha_{_{\G}}(y)) \rho_{_V}(x))\cr
&=&(\rho_{_U}(\alpha_{_{\G}}(x))(\rho_{_U}(y)))\otimes (\varphi_{_V}\circ  \varphi_{_V})+
(\varphi_{_U}\circ \varphi_{_U})\otimes (\rho_{_V}(\alpha_{_{\G}}(x))\rho_{_V}(y))\cr
&-&(\rho_{_U}(\alpha_{_{\G}}(y))\rho_{_U}(x))\otimes (\varphi_{_V}\circ \varphi_{_V})-
(\varphi_{_U}\circ \varphi_{_U})\otimes (\rho_{_V}(\alpha_{_{\G}}(y)) \rho_{_V}(x))
=(\rho_{_U}(\alpha_{_{\G}}(x))(\rho_{_U}(y))\cr 
&-&(\rho_{_U}(\alpha_{_{\G}}(y))\rho_{_U}(x)))\otimes (\varphi_{_V}\circ  \varphi_{_V})
+(\varphi_{_U}\circ \varphi_{_U})\otimes (\rho_{_V}(\alpha_{_{\G}}(x))\rho_{_V}(y)-(\rho_{_V}(\alpha_{_{\G}}(y)) \rho_{_V}(x)))\cr
&=&(\rho_{_U}[x,y]\varphi_{_U})\otimes (\varphi_{_V}\circ  \varphi_{_V})+
(\varphi_{_U}\circ \varphi_{_U})\otimes (\rho_{_V}[x,y]\varphi_{_V})\cr
&=&(\rho_{_U}[x,y]\otimes \varphi_{_V}+\varphi_{_V}\otimes \rho_{_V}[x,y] ) \circ (\varphi_{_U}\circ  \varphi_{_V})
 = 
 (\rho_{_U}\otimes \varphi_{_V}+\varphi_{_U}\otimes \rho_{_V})[x,y] (\varphi_{_U}\circ  \varphi_{_V}).
\eeqs
From this, we deduce that the linear map $\rho_{_U}\otimes \varphi_{_V}+\varphi_{_U} \otimes\rho_{_V}$
satisfies the relations \eqref{eq_rep_Hom_Lie1} and \eqref{eq_rep_Hom_Lie2},  i.e.  
$\rho_{_U}\otimes \varphi_{_V}+\varphi_{_U} \otimes\rho_{_V}$ is a representation of hom-Lie algebra 
$(\G, [,]_{_{\G}}, \alpha_{_{\G}})$.
$\cqfd$

For any linear map $\phi: U\rightarrow V,$ where $U$ and $V$ are  finite dimensional vector spaces,  we denote
by $\phi^*$ its dual  map defined by:
\beqs
\phi^*: 
\begin{array}{lll}
V^*\rightarrow U^*\\
v^* \mapsto \phi^*(v^*):
\begin{array}{lll}
U\rightarrow \K \\
u\mapsto <u, \phi^*(v^*)>=<\phi(u), v^*>,
\end{array}
\end{array} 
\eeqs
where $<,>$ is a natural pairing between 
$U$ and its dual space $U^*$. 
\begin{definition}
Consider a hom-Lie algebra $(\G, [, ], \alpha)$ and  a representation $(\rho, V)$ of $\G$.
A $1$-hom-cocycle $\delta: \G\rightarrow \G $ associated to  the linear map $\rho:\G \rightarrow V$ satisfies  the following relation:
\beq
\delta(\alpha([x,y]))=\rho(x)\delta(y)-\rho(y)\delta(x), \forall x, y\in \G; \qquad \forall x, y\in \G.
\eeq  
\end{definition}
\begin{theorem}\label{Theo_Hom_lll}
Let $(\A, \cdot, \alpha)$ be a hom-center-symmetric algebra  given by the  product 
$\beta^*: \A\otimes\A\rightarrow \A$
such that $\alpha^2=\id$. Suppose  there is another hom-center-symmetric algebra structure $"\circ "$ 
on its dual space $\A^*$ given by a linear map $\gamma^*: \A\otimes\A^*\rightarrow \A^*$.
Then   $(\G(\A), \G(\A^*), -\ad_{\cdot}^*, -\ad_{\circ}^*, \alpha, \alpha^*)$ is a matched pair of 
hom-Lie algebras $(\G(\A), \alpha)$ and $(\G(\A^*), \alpha^*)$ if and only if 
$\gamma: \A\rightarrow \A\otimes \A$ is a 1-hom-cocycle of the hom-Lie algebra  $(\G(\A), \alpha)$
associated to  $-(\ad_{\cdot}\otimes \alpha+\alpha\otimes\ad_{\cdot})$ and 
$\beta:\A^*\rightarrow \A^*\otimes\A^*$ is a 1-Hom-cocycle of $(\G(\A^*), \alpha^*)$ associated to
$-(\ad_{\circ}\otimes \alpha^*+\alpha^*\otimes\ad_{\circ})$.
\end{theorem}
\textbf{Proof}.

Let $\{e_1, e_2, \cdots, e_n\}$ be a basis of $\A,$ and $\{e_1^*, e_2^*, \cdots, e_n^*\}$  its dual basis.

For all $i,j\in \{1,2, \cdots,n \}$, consider  $\displaystyle e_i\cdot e_j=\sum_{k=1}^{n} c_{ij}^k e_k $ and 
$\displaystyle e_i^*\circ e_j^*=\sum_{k=0}^{n}f_{ij}^ke_k^*$, where $c_{ij}^k$ and  $f_{ij}^k\in \K$ are 
structure constants associated to the bilinear operations $\cdot $ and $\circ$, respectively. 
Then, for all $i, j, k\in \{1, \cdots, n\},$ we have:
\beqs
\left<\gamma(e_k), e_i^*\otimes e_j^*\right>=
\left<e_k, \gamma^*(e_i^*\otimes e_j^*)\right>=
\left<e_k, e_i^*\circ e_j^*\right>=
\left<e_k, \sum_{l=1}^{n}f_{ij}^l e_l^*\right>=
\sum_{l=1}^{n}f_{ij}^l\delta_l^k=
f_{ij}^k\\
=f_{ij}^k\left<e_i\otimes e_j, e_i^*\otimes e_j^*\right>=
\sum_{p=1}^{n}\sum_{q=1}^{n} f_{pq}^k\left<e_p, e_i^*\right>\left<e_q, e_j^*\right>=
\left<\sum_{p,q=1}^{n}f_{pq}^k e_p\otimes e_q, e_i^*\otimes e_j^*\right>,
\eeqs
and
\beqs
\left<\beta(e_k^*), e_i\otimes e_j\right>=
\left<e_k^*, \beta^*(e_i\otimes e_j)\right>=
\left< e_k^*, e_i \cdot e_j\right>=
\left<e_k^*, \sum_{l=1}^{n}c_{ij}^le_l\right>=
\sum_{l=1}^{n}c_{ij}^l \delta_l^k=c_{ij}^k \\
=c_{ij}^k\left<e_i\otimes e_j, e_i^*\otimes e_j^*\right>=
\sum_{p=1}^{n}\sum_{q=1}^{n} c_{pq}^k \left< e_p, e_i^*\right> \left<e_q, e_j^* \right>=
\left<\sum_{p,q=1}^{n}c_{pq}^k e_p^*\otimes e_q^*, e_i\otimes e_j\right>.
\eeqs
It follows that
$\displaystyle \gamma(e_k)=\sum_{i,j=1}^{n}f_{ij}^k e_i\otimes e_j$ and 
$\displaystyle \beta(e_k^*)=\sum_{i,j=1}^{n}c_{ij}^k e_i^*\otimes e_j^*$.
In addition, let $\displaystyle \alpha(e_i)=\sum_{k=1}^{n}d_i^ke_k$ and  
$\displaystyle \alpha^*(e_i^*)=\sum_{k=1}^{n}d^*{_i^ke_k^*}$. From the identity $\alpha^2=\id$,  we get 
$\displaystyle \sum_{k=1}^{n}\sum_{l=1}^{n}d_i^kd_k^l e_l=\sum_{l=1}^{n}\delta_i^le_l=e_i$, with
$d_i^kd_k^l=\delta_i^l$. Besides, we have: $<\alpha^*(e_i^*), e_j>=d_i^*{^j}=<e_i^*, \alpha(e_j)>=d_j^i$
which implies $d_i^*{^j}=d_j^i$.
 Furthermore, we also have 
\beq\label{eq_key1}
 \gamma(\alpha([e_i, e_j]))=\sum_{k,l,m,p=1}^{n}
 \{f_{mp}^l (c_{ij}^k-c_{ji}^k)d_k^l\}  e_m\otimes e_p,
\eeq
and 
\beq\label{eq_key2} 
 \beta(\alpha^*([e_i^*, e_j^*]))=\sum_{k,l,m,p=1}^{n} 
\{ c_{mp}^l (f_{ij}^k-f_{ji}^k)d^*{_k^l}\} e_m^*\otimes e_p^*.
\eeq
Besides, we have for all $i, j \in \{1, 2, \cdots, n \}:$
\beqs
&&\{(-\ad_{\cdot}{e_i})\otimes \alpha+\alpha\otimes(-\ad_{\cdot}{e_i})\}\gamma(e_j)-
\{(-\ad_{\cdot}{e_j})\otimes \alpha+\alpha\otimes(-\ad_{\cdot}{e_j})\}\gamma(e_i)\cr 
&&=\{(\ad_{\cdot}{e_j}\otimes \alpha)\gamma(e_i)-(\ad_{\cdot}{e_i}\otimes \alpha)\gamma(e_j)\}+
\{(\alpha\otimes\ad_{\cdot}{e_j})\gamma(e_i)-(\alpha\otimes\ad_{\cdot}{e_i})\gamma(e_j)\}
\cr
&&=\sum_{k, l=1}^{n}\{ f_{kl}^i(\ad_{\cdot}{e_j}\otimes \alpha)e_k\otimes e_l-
 f_{kl}^{j} (\ad_{\cdot}{e_i}\otimes \alpha)e_k\otimes e_l\} 
 \cr 
&&
+\sum_{l,k=1}^{n} \{ f_{lk}^i (\alpha\otimes\ad_{\cdot}{e_j})e_l\otimes e_k-
f_{lk}^j(\alpha\otimes\ad_{\cdot}{e_i})e_l\otimes e_k\} 
\cr
&&
=\sum_{k, l=1}^{n}\{ f_{kl}^i([e_j,e_k]\otimes \alpha(e_l))-
 f_{kl}^{j} ([e_i, e_k]\otimes \alpha(e_l))\}
  \cr 
&&
+\sum_{l,k=1}^{n} \{ f_{lk}^i (\alpha(e_l)\otimes[e_j, e_k]-
f_{lk}^j(\alpha (e_l)\otimes [e_i,e_k]\} 
\cr
&&=\sum_{k, l=1}^{n}\left\lbrace 
\sum_{m=1}^{n}
\{
 f_{kl}^i(c_{jk}^m-c_{kj}^m)e_m\otimes \alpha(e_l)-
 f_{kl}^{j}(c_{ik}^m-c_{ki}^m)e_m\otimes \alpha(e_l)\}
 \right\rbrace\cr  
&&+\sum_{l,k=1}^{n} \left\lbrace
\sum_{p=1}^{n} \{
f_{lk}^i (c_{jk}^p-c_{kj}^p)\alpha(e_l)\otimes e_p-
f_{lk}^j(c_{ik}^p-c_{ki}^p)\alpha (e_l)\otimes e_p
\}
\right\rbrace 
\cr
&&=\sum_{k, l=1}^{n} 
\sum_{m=1}^{n}
\left\lbrace
\sum_{p=1}^{n}
\{
 f_{kl}^i(c_{jk}^m-c_{kj}^m)-
 f_{kl}^{j} (c_{ik}^m-c_{ki}^m)\}d_l^p e_m\otimes e_p 
 \right\rbrace\cr  
&&+\sum_{l,k=1}^{n} 
\sum_{p=1}^{n} 
\left\lbrace
\sum_{m=1}^{n}
\{
f_{lk}^i (c_{jk}^p-c_{kj}^p)-
f_{lk}^j(c_{ik}^p-c_{ki}^p)\} d_l^m  e_m \otimes e_p
\right\rbrace 
\cr
&&=\sum_{k, l, m, p =1}^{n} 
\left\lbrace
\{
 f_{kl}^i(c_{jk}^m-c_{kj}^m)-
 f_{kl}^{j} (c_{ik}^m-c_{ki}^m)\}d_l^p  
+ 
\{
f_{lk}^i (c_{jk}^p-c_{kj}^p)-
f_{lk}^j(c_{ik}^p-c_{ki}^p)\} d_l^m  
 \right\rbrace
e_m
 \otimes e_p
\eeqs
Therefore, 
\beqs
\{(-\ad_{\cdot}{e_i})\otimes \alpha+\alpha\otimes(-\ad_{\cdot}{e_i})\}\gamma(e_j)-
\{(-\ad_{\cdot}{e_j})\otimes \alpha+\alpha\otimes(-\ad_{\cdot}{e_j})\}\gamma(e_i)
\cr
=\sum_{k, l, m, p =1}^{n} 
\left\lbrace
\{
 f_{kl}^i(c_{jk}^m-c_{kj}^m)-
 f_{kl}^{j} (c_{ik}^m-c_{ki}^m)\}d_l^p  
+  
\{
f_{lk}^i (c_{jk}^p-c_{kj}^p)-
f_{lk}^j(c_{ik}^p-c_{ki}^p)\} d_l^m  
 \right\rbrace
e_m
 \otimes e_p
\eeqs
By using the fact that $\gamma$ is a $1$-cocycle associated to the underlying hom-Lie algebra
$(\G(\A), \alpha)$ with the representation
$-(\ad_{\cdot}\otimes \alpha+\alpha\otimes\ad_{\cdot}), $
and taking into account the relation \eqref{eq_key1}, we get:
\beq\label{eq_hom_long_1}
\sum_{k,l,m,p=1}^{n}
 f_{mp}^l (c_{ij}^k-c_{ji}^k)d_k^l&=&
 \sum_{k, l, m, p =1}^{n} 
 \{
  (f_{kl}^i(c_{jk}^m-c_{kj}^m)-
  f_{kl}^{j} (c_{ik}^m-c_{ki}^m))d_l^p  
 \\ &+& 
( f_{lk}^i (c_{jk}^p-c_{kj}^p)
- f_{lk}^j(c_{ik}^p-c_{ki}^p))
 \} d_l^m  
 \}\nonumber
\eeq
Similarly, we  obtain
\beqs
&&\{(-\ad_{\circ}^*{e_i^*})\otimes \alpha^*+\alpha^*\otimes(-\ad_{\circ}^*{e_i^*})\}\gamma(e_j^*)-
 \{(-\ad_{\circ}^*{e_j^*})\otimes \alpha^*+\alpha^*\otimes(-\ad_{\circ}^*{e_j^*})\}\gamma(e_i^*)\cr 
&=& \{(\ad_{\circ}^*{e_j^*}\otimes \alpha^*)\gamma(e_i^*)-(\ad_{\circ}^*{e_i^*}\otimes \alpha^*)\gamma(e_j^*)\}+
 \{(\alpha^*\otimes\ad_{\circ}^*{e_j^*})\gamma(e_i^*)-(\alpha^*\otimes\ad_{\circ}^*{e_i^*})\gamma(e_j^*)\}
 \cr
&=& \sum_{k, l=1}^{n}\{ c_{kl}^i(\ad_{\circ}{e_j^*}\otimes \alpha^*)e_k^*\otimes e_l^*-
  c_{kl}^{j} (\ad_{\circ}{e_i^*}\otimes \alpha^*)e_k^*\otimes e_l^*\}\cr 
&+& \sum_{l,k=1}^{n} \{ c_{lk}^i (\alpha^*\otimes\ad_{\circ}{e_j^*})e_l^*\otimes e_k^*-
 c_{lk}^j(\alpha^*\otimes\ad_{\circ}{e_i^*})e_l^*\otimes e_k^*\}  \cr
&=& \sum_{k, l=1}^{n}\{ c_{kl}^i([e_j^*,e_k^*]\otimes \alpha^*(e_l^*))-
  c_{kl}^{j} ([e_i^*, e_k^*]\otimes \alpha^*(e_l^*))\} \cr 
&+&\sum_{l,k=1}^{n} \{ f_{lk}^i (\alpha^*(e_l)\otimes[e_j, e_k]-
 c_{lk}^j(\alpha^* (e_l^*)\otimes [e_i^*,e_k^*]\} 
 \cr
&=& \sum_{k, l=1}^{n}\left\lbrace 
 \sum_{m=1}^{n}
 \{
  c_{kl}^i(f_{jk}^m-f_{kj}^m)e_m^*\otimes \alpha(e_l^*)-
  c_{kl}^{j}(f_{ik}^m-f_{ki}^m)e_m^*\otimes \alpha(e_l^*)\}
  \right\rbrace\cr  
&+&\sum_{l,k=1}^{n} \left\lbrace
 \sum_{p=1}^{n} \{
 c_{lk}^i (f_{jk}^p-f_{kj}^p)\alpha(e_l^*)\otimes e_p^*-
 c_{lk}^j(f_{ik}^p-f_{ki}^p)\alpha (e_l^*)\otimes e_p^*
 \}
 \right\rbrace 
 \cr
&=& \sum_{k, l=1}^{n} 
 \sum_{m=1}^{n}
 \left\lbrace
 \sum_{p=1}^{n}
 \{
  c_{kl}^i(f_{jk}^m-f_{kj}^m)-
  c_{kl}^{j} (f_{ik}^m-f_{ki}^m)\}d^*{_l^p} e_m^*\otimes e_p^* 
  \right\rbrace\cr  
&+& \sum_{l,k=1}^{n} 
 \sum_{p=1}^{n} 
 \left\lbrace
 \sum_{m=1}^{n}
 \{
 c_{lk}^i (f_{jk}^p-f_{kj}^p)-
 c_{lk}^j(f_{ik}^p-f_{ki}^p)\} d^*{_l^m}  e_m^* \otimes e_p^*
 \right\rbrace 
 \cr
&=& \sum_{k, l, m, p =1}^{n} 
 (
 \{
  c_{kl}^i(f_{jk}^m-f_{kj}^m)-
  c_{kl}^{j} (f_{ik}^m-f_{ki}^m)\}d^*{_l^p}  
 +  
 \{
  c_{lk}^i (f_{jk}^p-f_{kj}^p)-
 c_{lk}^j(f_{ik}^p-f_{ki}^p)\} d^*{_l^m}  
  )
 e_m^*
  \otimes e_p^*
\eeqs
Therefore, 
\beqs
 \{(-\ad^*_{\circ}{e_i^*})\otimes \alpha^*+\alpha^*\otimes(-\ad^*_{\circ}{e_i^*})\}\gamma(e_j^*)-
 \{(-\ad^*_{\circ}{e_j^*})\otimes \alpha^*+\alpha^*\otimes(-\ad^*_{\circ}{e_j^*})\}\gamma(e_i^*)=
 \cr
 \sum_{k, l, m, p =1}^{n} 
 (
 \{
  c_{kl}^i(f_{jk}^m-f_{kj}^m)-
  c_{kl}^{j} (f_{ik}^m-f_{ki}^m)\}d^*{_l^p}  
 +   
 \{
 c_{lk}^i (f_{jk}^p-f_{kj}^p)-
 c_{lk}^j(f_{ik}^p-f_{ki}^p)\} d^*{_l^m}  
  )
 e_m^*
  \otimes e_p^*
\eeqs
By using the fact that $\gamma$ is a $1$-cocycle associated to the underlying hom-Lie algebra
$(\G(\A^*), \alpha^*)$  with the representation
$-(\ad_{\circ}\otimes \alpha^*+\alpha^*\otimes\ad_{\circ}), $
and using the relation \eqref{eq_key2}, we have:
\beq\label{eq_hom_long_2}
\sum_{k,l,m,p=1}^{n}
 c_{mp}^l (f_{ij}^k-f_{ji}^k)d^*{_k^l}&=&
 \sum_{k, l, m, p =1}^{n} 
 \{
  (c_{kl}^i(f_{jk}^m-f_{kj}^m)-
  c_{kl}^{j} (f_{ik}^m-f_{ki}^m))d^*{_l^p}\cr
   &+&  
( c_{lk}^i (f_{jk}^p-f_{kj}^p) 
 -c_{lk}^j(f_{ik}^p-f_{ki}^p)) d^*{_l^m}
 \}.\nonumber
\eeq
Besides,  for all $i, j, k\in \{1, 2, \cdots , n\},$ we obtain:
\beqs
\left<\ad_{\cdot}^*(\alpha(e_i))e_j^*, e_k
\right>=
\left<e_j^*, [\alpha(e_i), e_k]
\right>
&=&\sum_{l=1}^{n}d_i^l
\left<e_j^*, [e_l, e_k]
\right>
=\sum_{l,p=1}^{n}d_i^l(c_{lk}^p-c_{kl}^p)
\left<e_j^*, e_p
\right>\cr 
&=& 
\sum_{l=1}^{n}d_i^l(c_{lk}^j-c_{kl}^j)
=
\left<
\sum_{l,p=1}^{n}d_i^l(c_{lp}^j-c_{pl}^j)e_p^*, e_k
\right>,
\eeqs
and then we get
\beqs
\ad_{\cdot}^*(\alpha(e_i))e_j^*=\sum_{k,l=1}^{n}d_i^k(c_{kl}^j-c_{lk}^j)e_l^*,\qquad 
\ad_{\circ}^*(\alpha^*(e_i^*))e_j=\sum_{k,l=1}^{n}d{^*}_i^k(f_{kl}^j-f_{lk}^j)e_l,
\eeqs
\beq
\ad_{\cdot}^*(e_i)e_j^*=\sum_{k=1}^{n}(c_{ik}^j-c_{ki}^j)e_k^*,
\eeq
and  
\beq
\ad_{\circ}^*(e_i^*)e_j=\sum_{k=1}^{n}(f_{ik}^j-f_{ki}^j)e_k.
\eeq
In addition, we also have
\beq
\ad_{\circ}^*(\alpha^*(e_i^*))[e_k, e_j] &=&\sum_{l,m,p=1}^{n} d{^*}_i^p(c_{kj}^l-c_{jk}^l)(f_{pm}^l-f_{mp}^l) e_m \cr
\ad_{\circ}^*(\alpha^*(e_i^*))[e_k, e_j] &=&\sum_{l,m,p=1}^{n} d_p^i(c_{kj}^l-c_{jk}^l)(f_{pm}^l-f_{mp}^l) e_m
\eeq
\beq
\ad_{\cdot}^*(\alpha(e_i))[ e_k^*,e_j^*]
&=& \sum_{l,m,p=1}^{n}d_i^p(f_{kj}^l-f_{jk}^l)(c_{pm}^l-c_{mp}^l) e_m^*\cr 
\ad_{\cdot}^*(\alpha(e_i))[ e_k^*,e_j^*]&=&
\sum_{l,m,p=1}^{n}d_i^*{^p}(f_{kj}^l-f_{jk}^l)(c_{pm}^l-c_{mp}^l) e_m^*
\eeq
Then, we find:
\beqs
&&\ad_{\circ}^*(\ad_{\cdot}^*(e_k)e_i^*)\alpha(e_j)-
\ad_{\circ}^*(\ad_{\cdot}^*(e_j)e_i^*)\alpha(e_k)-
[\ad_{\circ}^*(e_i^*)e_j, \alpha(e_k)]-
[\alpha(e_j), \ad_{\circ}^*(e_i^*)e_k]\cr 
&=&\ad_{\circ}^*(\ad_{\cdot}^*(e_k)e_i^*)\alpha(e_j)-
\ad_{\circ}^*(\ad_{\cdot}^*(e_j)e_i^*)\alpha(e_k)+
\alpha(e_k)\cdot (\ad_{\circ}^*(e_i^*)e_j)-
(\ad_{\circ}^*(e_i^*)e_j)\cdot\alpha(e_k)\cr 
&+& (\ad_{\circ}^*(e_i^*)e_k)\cdot\alpha(e_j)-
\alpha(e_j)\cdot (\ad_{\circ}^*(e_i^*)e_k)\cr &=& 
\{
\ad_{\circ}^*(\ad_{\cdot}^*(e_k)e_i^*)\alpha(e_j)-
\alpha(e_j)\cdot (\ad_{\circ}^*(e_i^*)e_k)+
(\ad_{\circ}^*(e_i^*)e_k)\cdot\alpha(e_j)
\}\cr 
&+&\{
\alpha(e_k)\cdot (\ad_{\circ}^*(e_i^*)e_j)
\ad_{\circ}^*(\ad_{\cdot}^*(e_j)e_i^*)\alpha(e_k)-
(\ad_{\circ}^*(e_i^*)e_j)\cdot\alpha(e_k)
\}\cr&=& 
\sum_{l=1}^{n}
\{
d_j^l\{\ad_{\circ}^*(\ad_{\cdot}^*(e_k)e_i^*)e_l-
e_l\cdot(\ad_{\circ}^*(e_i^*)e_k)+
(\ad_{\circ}^*(e_i^*)e_k)\cdot e_l
\}\cr 
&+& d_k^l\{e_l\cdot(\ad_{\circ}^*(e_i^*)e_j)-\ad_{\circ}^*(\ad_{\cdot}^*(e_j)e_i^*)e_l-
(\ad_{\circ}^*(e_i^*)e_j)\cdot e_l 
\}
\}\cr &=&
\sum_{l, r=1}^{n}\{
d_j^l\{(c_{kr}^i-c_{rk}^i)\ad_{\circ}^*(e_r^*)e_l-
(f_{ir}^k-f_{ri}^k)(e_l\cdot e_r)+(f_{ir}^k-f_{ri}^k)(e_r\cdot e_l)
\}
 \cr &+& 
d_k^l\{(f_{ir}^j-f_{ri}^j)(e_l\cdot e_r)-(c_{jr}^i-c_{rj}^i)\ad_{\circ}^*(e_r^*)e_l-
(f_{ir}^j-f_{ri}^j)(e_r\cdot e_j)
\}
\}\cr&=& 
\sum_{l,r,m=1}^{n}\{
d_j^l(c_{kr}^i-c_{rk}^i)(f_{rm}^l-f_{mr}^l)+
d_j^l(f_{ir}^k-f_{ri}^k)(c_{rl}^m-c_{lr}^m)+
d_k^l(f_{ir}^j-f_{ri}^j)(c_{lr}^m-c_{rl}^m)\cr 
&+& d_k^l(c_{rj}^i-c_{jr}^i)(f_{rm}^l-f_{mr}^l)
\}e_m,
\eeqs
and the following equalities hold:
\beqs
&& \sum_{l,m,p=1}^{n} d_p^i(c_{kj}^l-c_{jk}^l)(f_{pm}^l-f_{mp}^l)=
 \sum_{l,r,m=1}^{n}
 \{
 d_j^l(c_{kr}^i-c_{rk}^i)(f_{rm}^l-f_{mr}^l)\cr &+&
 d_j^l(f_{ir}^k-f_{ri}^k)(c_{rl}^m-c_{lr}^m)+ d_k^l(f_{ir}^j-f_{ri}^j)(c_{lr}^m-c_{rl}^m)+ 
 d_k^l(c_{rj}^i-c_{jr}^i)(f_{rm}^l-f_{mr}^l)
 \}\Longleftrightarrow \cr 
 &&
 \sum_{l,m,p=1}^{n} d_p^i(c_{kj}^l-c_{jk}^l)f_{pm}^l-
 \sum_{l,m,p=1}^{n} d_p^i(c_{kj}^l-c_{jk}^l)f_{mp}^l\cr 
 &=&
 \sum_{l,r,m=1}^{n}\{
 d_j^l(c_{kr}^i-c_{rk}^i)f_{rm}^l +
 d_j^lf_{ir}^k(c_{rl}^m-c_{lr}^m) 
 +d_k^lf_{ir}^j(c_{lr}^m-c_{rl}^m)+ 
 d_k^l(c_{rj}^i-c_{jr}^i)f_{rm}^l \cr &-&
 \sum_{l,r,m=1}^{n}\{
 d_j^l(c_{kr}^i-c_{rk}^i)f_{mr}^l+
 d_j^lf_{ri}^k(c_{rl}^m-c_{lr}^m) 
 +d_k^lf_{ri}^j(c_{lr}^m-c_{rl}^m)+  
 d_k^l(c_{rj}^i-c_{jr}^i)f_{mr}^l\} \Longleftrightarrow \cr 
 &&
 \sum_{l,m,p=1}^{n} d_p^i(c_{kj}^l-c_{jk}^l)f_{pm}^l 
 - \sum_{l,r,m=1}^{n}\{
 d_j^l(c_{kr}^i-c_{rk}^i)f_{rm}^l+ d_j^lf_{ir}^k(c_{rl}^m-c_{lr}^m)+ 
 d_k^lf_{ir}^j(c_{lr}^m-c_{rl}^m)\cr &+& 
 d_k^l(c_{rj}^i-c_{jr}^i)f_{rm}^l\}= 
 \sum_{l,m,p=1}^{n} d_p^i(c_{kj}^l-c_{jk}^l)f_{mp}^l- \sum_{l,r,m=1}^{n}\{
 d_j^l(c_{kr}^i-c_{rk}^i)f_{mr}^l+
 d_j^lf_{ri}^k(c_{rl}^m-c_{lr}^m)\cr &+&  
 d_k^lf_{ri}^j(c_{lr}^m-c_{rl}^m) +  
 d_k^l(c_{rj}^i-c_{jr}^i)f_{mr}^l\} \Longleftrightarrow 
 \sum_{l,m,p=1}^{n} d_p^i(c_{kj}^l-c_{jk}^l)f_{pm}^l=
 \sum_{l,r,m=1}^{n} 
 d_j^l(c_{kr}^i-c_{rk}^i)f_{rm}^l\cr &+& 
 d_j^lf_{ir}^k(c_{rl}^m-c_{lr}^m) + 
 d_k^lf_{ir}^j(c_{lr}^m-c_{rl}^m) 
 + d_k^l(c_{rj}^i-c_{jr}^i)f_{rm}^l. 
\eeqs
Therefore, we obtain
\beq\label{eq_Hom_center_ll1}
\sum_{l,m,p=1}^{n} d_p^i(c_{kj}^l-c_{jk}^l)f_{pm}^l&=&
\sum_{l,r,m=1}^{n}
d_j^l(c_{kr}^i-c_{rk}^i)f_{rm}^l+  
d_j^lf_{ir}^k(c_{rl}^m-c_{lr}^m)\cr 
&+& 
d_k^lf_{ir}^j(c_{lr}^m-c_{rl}^m) +d_k^l(c_{rj}^i-c_{jr}^i)f_{rm}^l. 
\nonumber
\eeq
Similarly, we compute:
\beqs
&&\ad_{\cdot}^*(\ad_{\circ}^*(e_k^*)e_i)\alpha^*(e_j^*)-
\ad_{\cdot}^*(\ad_{\circ}^*(e_j^*)e_i)\alpha^*(e_k^*)-
[\ad_{\cdot}^*(e_i)e_j^*, \alpha^*(e_k^*)]-
[\alpha^*(e_j^*), \ad_{\cdot}^*(e_i)e_k^*]\cr 
&=& 
\sum_{l,r,m=1}^{n}
d_j^*{^l}(f_{kr}^i-f_{rk}^i)(c_{rm}^l-c_{mr}^l)+
d_j^*{^l}(c_{ir}^k-c_{ri}^k)(f_{rl}^m-f_{lr}^m)+
d_k^*{^l}(c_{ir}^j-c_{ri}^j)(f_{lr}^m-f_{rl}^m) \cr 
&+& d_k^*{^l}(f_{rj}^i-f_{jr}^i)(c_{rm}^l-c_{mr}^l)
\Longleftrightarrow
\sum_{l,m,p=1}^{n}d_i^*{^p}(f_{kj}^l-f_{jk}^l)c_{pm}^l-
\sum_{l,m,p=1}^{n}d_i^*{^p}(f_{kj}^l-f_{jk}^l)c_{mp}^l
\cr &=&  
\sum_{l,r,m=1}^{n}
d_j^*{^l}(f_{kr}^i-f_{rk}^i)c_{rm}^l+ 
d_j^*{^l}c_{ir}^k(f_{rl}^m-f_{lr}^m)  
+d_k^*{^l}c_{ir}^j(f_{lr}^m-f_{rl}^m)+    
d_k^*{^l}(f_{rj}^i-f_{jr}^i)c_{rm}^l\cr  
&-&\sum_{l,r,m=1}^{n}\{
d_j^*{^l}(f_{kr}^i-f_{rk}^i)c_{mr}^l+   
d_j^*{^l}c_{ri}^k(f_{rl}^m-f_{lr}^m)  
+d_k^*{^l}c_{ri}^j(f_{lr}^m-f_{rl}^m)+    
d_k^*{^l}(f_{rj}^i-f_{jr}^i)c_{mr}^l \}\Longleftrightarrow\cr 
&&\sum_{l,m,p=1}^{n}d_i^*{^p}(f_{kj}^l-f_{jk}^l)c_{pm}^l=
\sum_{l,r,m=1}^{n}
d_j^*{^l}(f_{kr}^i-f_{rk}^i)c_{rm}^l+ 
d_j^*{^l}c_{ir}^k(f_{rl}^m-f_{lr}^m)  
+d_k^*{^l}c_{ir}^j(f_{lr}^m-f_{rl}^m)   \cr
&+& d_k^*{^l}(f_{rj}^i-f_{jr}^i)c_{rm}^l
\eeqs
Therefore, we get
\beq\label{eq_Hom_center_ll2}
\sum_{l,m,p=1}^{n}d_i^*{^p}(f_{kj}^l-f_{jk}^l)c_{pm}^l&=&
\sum_{l,r,m=1}^{n}
d_j^*{^l}(f_{kr}^i-f_{rk}^i)c_{rm}^l+ 
d_j^*{^l}c_{ir}^k(f_{rl}^m-f_{lr}^m) \\
&+& d_k^*{^l}c_{ir}^j(f_{lr}^m-f_{rl}^m)  
+d_k^*{^l}(f_{rj}^i-f_{jr}^i)c_{rm}^l.\nonumber 
\eeq
By adding the relation $\alpha^2=\id$, we obtain that 
\eqref{eq_hom_long_1} is equivalent to \eqref{eq_Hom_center_ll1}, and  using 
$\alpha^*{^2}=\id$,  the equation
\eqref{eq_hom_long_1} is equivalent to the equation  \eqref{eq_Hom_center_ll2}.
$\cqfd$
\begin{definition}
Let $(\A, \alpha)$ be a vector space. A hom-center-symmetric bialgebra structure  on $\A$ is a pair
of linear map $(\gamma, \beta)$ with $\alpha\in \mathfrak{gl}(\A)$, $\gamma: \A\rightarrow \A\otimes\A$, 
$\beta: \A^*\rightarrow \A^*\otimes\A^*$ and
\begin{itemize}
\item $\gamma^*: \A^*\otimes\A^*\rightarrow \A^*$ is a hom-center-symmetric algebra structure algebra on $\A^*$,
\item $\beta^*: \A\otimes\A\rightarrow \A$ is a hom-center symmetric algebra structure algebra on $\A$, 
\item $\gamma$ is a 1-hom-cocycle of $\G(\A)$ associated to  $-(\ad_{\cdot}\otimes \alpha+\alpha\otimes\ad_{\cdot})$
 with values in $\A \otimes\A$,
\item $\beta$ is a 1-hom-cocycle of $\G(\A^*)$ associated to  
$-(\ad_{\circ}\otimes \alpha^*+\alpha^*\otimes\ad_{\circ})$
 with values in $\A^*\otimes\A^*$.
\end{itemize}
We denote this hom-center-symmetric bialgebra by $(\A, \A^*, \gamma, \beta, \alpha)$.
\end{definition}
\begin{definition} 
A Manin triple of hom-center-symmetric algebras 
$(\A,\cdot , \alpha)$ and $(\B, \ast, \alpha')$ is a triple 
$(\A\oplus\B, \A, \B),$ together with a nondegenerate symmetric 
bilinear form $\bb(\cdot, \cdot )$ on the hom-center-symmetric algebra 
$(\A\oplus\B, \star, \alpha\oplus\alpha')$ such that 
\begin{itemize}
\item
$\bb$ is invariant, i. e.  for all $x,y,z\in \A$ and $a, b,c \in \B$,
\beqs
\bb((x+a)\star (y+b), (z+c))&=&\bb((x+a),(y+b)\star (z+c)),\cr 
\bb((\alpha\oplus\alpha')(x+a), y+b)&=&\bb(x+a, (\alpha\oplus\alpha')(y+b)).
\eeqs
\item The hom-center-symmetric algebras $\A$ and $\B$ are isotropic hom-center-symmetric algebras of  $A\oplus\B$.
\end{itemize}
\end{definition}
In particular, if $(\A, \cdot, \alpha)$ is a hom-center-symmetric algebra, and if there exists a hom-center-symmetric 
algebra structure on its dual space $\A^*$ denoted  by $(\A^*,\circ, \alpha^*)$, then there is a hom-center-symmetric 
algebra on the direct sum of the underlying vector space of $\A$ and its dual space $\A^*$ such that
$(\A\oplus\A^*, \A, \A^*)$ is the associated Manin triple with the invariant bilinear symmetric form given by
$\bb_{\A}(x+a, y+b)=<x,b>+<y, a>$ for all $x, y\in \A$ and $a, b\in \A^*.$ It is
called the standard Manin triple of the hom-center-symmetric algebra $\A;$  $<,>$ 
is a natural paring between algebra and its dual space.
\begin{proposition}\label{Prop_Hom_dd}
Let $(\A, \cdot, \alpha)$ and $(\A^*, \circ, \alpha^*)$ be two hom-center-symmetric algebras. Then,
$(\A, \A^*, R_{\cdot}^*, L_{\cdot}^*, R_{\circ}^*, L_{\circ}^*, \alpha, \alpha^*)$
 is a matched pair of the  hom-center-symmetric algebras $(\A, \cdot, \alpha)$ and $(\A^*, \circ, \alpha^*)$ 
 if and only if
$(\A\oplus\A^*, \A, \A^*)$ is a standard Manin triple.
\end{proposition}
\textbf{Proof}.

Let us compute and compare the following relations: $\displaystyle \bb_{\A}\left((x+a)\ast(y+b), (z+c)\right)$ 
and\\ $ \displaystyle
 \bb_{\A}\left((x+a), (y+b)\ast(z+c)\right)$  $\forall x, y, z \in \A$ and  $\forall a, b, c \in \A^{*}.$
\beqs
\bb_{\A} ((x+a)\ast(y+b), (z+c) ) 
&=&\bb_{\A}(xy+R_{\circ}^*(a)y+L_{\circ}^*(b)x+ 
a \circ b+R_{\cdot}^*(x)b+L_{\cdot}^*(y)a, z+c)\cr
&=& \left< xy+R^*_{\circ}(a)y+L_{\circ}^*(b)x, c  \right>+
\left< z, a \circ b + R^*_{\cdot}(x)b+L^*_{\cdot}(y)a  \right> \cr
&=& \left<xy, c  \right>+ \left<R^*_{\circ}(a)y, c \right>+ 
\left< L^*_{\circ}(b)x , c\right> 
+\left<z, a \circ b  \right>
+\left<z, R^*_{\cdot}(x)b \right>\cr 
&+&
\left<z, L^*_{\cdot}(y)a   \right> 
= 
\left<  xy, c \right>+
\left<y, R_a(c) \right>+
\left<x, L_b(c)  \right>+ 
\left< z, a \circ b  \right>\\
&+&
\left< R_x(z) , b \right>+
\left< L_y(z) , a \right>\\
& =&
\left< xy, c\right>+
\left<  y, c\circ a \right>+
\left< x , b \circ c  \right>+ 
\left< z, a \circ b  \right>+
\left<zx, b   \right>+
\left<yz, a \right>.
\eeqs
\beqs
 \bb_{\A}\left((x+a), (y+b)\ast(z+c)\right)
 &=& \bb_{\A}(x+a, 
yz+R_{\circ}^*(b)z+L^*_{\circ}(c)y+  b \circ c+
R^*_{\cdot}(y)c+L^*_{\cdot}(z)b)\cr
&=& \left<x, b \circ c +R^*_{\cdot}(y)c+L^*_{\cdot}(z)b \right>+ 
\left<  yz+R^*_{\circ}(b)z+L_{\circ}^*(c)y   ,a   \right> \cr
&+&
\left<x, b \circ c \right>+
\left<x, R^*_{\cdot}(y)c  \right>+
\left<x , L^*_{\cdot}(z)b  \right> 
+\left< yz , a\right>+ 
\left<R^*_{\circ}(b)z, a   \right>\\
&+&
\left<L^*_{\circ}(c)y, a  \right> =
\left<x, b\circ c   \right>+
\left<R_y(x),  c \right>+
\left<L_z(x), b  \right> \cr
&+&
\left< yz, a \right>+
\left<z, R_b(a) \right>+
\left<y, L_c(a) \right> \cr
&=&
\left<x, b\circ c  \right>+
\left<xy, c  \right>+
\left<zx, b  \right>+ 
\left<yz, a  \right>+
\left<z, a \circ b \right>+
\left<y, c \circ a  \right>.
\eeqs
It follows that 
\beq
\bb_{\A} ((x+a)\ast(y+b), (z+c) ) =\bb_{\A}\left((x+a), (y+b)\ast(z+c)\right)
\eeq 
which expresses the invariance of the standard bilinear form on $\A \oplus \A^*.$ 
Therefore, $(A\oplus \A^*, \A, \A^*)$ is the standard Manin triple of the center-symmetric algebras $\A$ and  $\A^*.$
$\cqfd$
\begin{theorem}
Let $(\A, \cdot , \alpha)$ be a hom-center-symmetric algebra and let $(\A^*, \circ, \alpha^*)$ be a 
hom-center-symmetric algebra structure on its dual space $\A^*$. Then the following conditions 
are equivalent:
\begin{enumerate}
\item[(i)]  $(\A \oplus \A^* , \A, \A^* )$ is the standard Manin triple of the hom-center-symmetric algebras 
$(\A, \cdot , \alpha)$ and $(\A^*, \circ, \alpha^*)$;
\item[(ii)]
$(\A, \A^* , R_{\cdot}^* , L_{\cdot}^* , R_{\circ}^* , L_{\circ}^*, \alpha, \alpha^* )$ is a matched pair of the
hom-center-symmetric algebras $(\A, \cdot , \alpha)$ and $(\A^*, \circ, \alpha^*)$;
\item[(iii)]
$(\G(\A), \G(\A^* ),-\ad_{\cdot}^*,-\ad_{\circ}^*,\alpha,\alpha^*)$ 
is a matched pair of the sub-adjacent Lie algebras $(\G(\A), \cdot , \alpha)$ and $(\G(\A^*), \circ, \alpha^*)$;
\item[(iv)]
$(\A, \A^*) $ is a center -symmetric bialgebra.
\end{enumerate}
\end{theorem}
\textbf{Proof}.

By the Proposition~\ref{Prop_Hom_dd}, $(i)\Longleftrightarrow (ii)$. From Theorem~\ref{Theo_Hom_ll}, 
 $(ii)\Longleftrightarrow (iii)$. According to  Theorem~\ref{Theo_Hom_lll},  $(iii)\Longleftrightarrow (iv)$. 
$\cqfd$
\begin{definition}
Let $(\A, \A^*, \gamma_{_{\A}},\beta_{_{\A}}, \alpha_{1})$ and 
$(\B, \B^*,\gamma_{_{\B}}, \beta_{_{\B}},  \alpha_{2})$ be two 
hom-center-symmetric bialgebras. A homomorphism of a  hom-center-symmetric bialgebra 
$f:\A\rightarrow \B$ is a homomorphism
of a hom-center-symmetric algebra such that $f^*:\B^*\rightarrow \A^*$ is also a homomorphism of a
hom-center-symmetric algebra, that is for all $x\in \A$, and  $a^*\in \B^*$, $f$ satisfies
\beq
(f\otimes f)\gamma_{_{\A}}(x)=\gamma_{_{\B}}(f(x)), \quad f\alpha_1=\alpha_2f, \quad \alpha_2^*f=f\alpha_1^*, \quad
(f^*\otimes f^*)\beta_{_{\B}}(a^*)=\beta_{_{\A}}(f^*(a^*)).
\eeq
\end{definition}
 \section{Concluding remarks}
In this work, we have constructed the hom-center-symmetric algebras,  and discussed
 their main  properties.  We have defined their 
bimodules and matched pairs. Furthermore, 
we  have established the  dual bimodules 
of  bimodules of   hom-center-symmetric algebras, and linked them  to the matched pairs 
of hom-Lie algebras. Besides, we have derived the Manin triple of 
hom-center-symmetric algebras. Finally, we have  provided    a theorem linking  hom-center-symmetric algebras,  associated matched pairs, 
hom-center-symmetric bialgebras, and  matched pairs of 
sub-adjacent hom-Lie algebras.

\end{document}